\documentclass[11pt]{amsart} 
\usepackage{times}
\usepackage{amsfonts,amssymb,amscd}
\DeclareMathOperator{\Area}{Area}
\DeclareMathOperator{\dist}{dist}
\DeclareMathOperator{\len}{Length}
\newtheorem{thm}{Theorem}
\newtheorem{defn}{Definition}
\newtheorem{lem}{Lemma}
\newtheorem{cor}{Corollary}
\newtheorem{prop}{Proposition}
\newtheorem{ram}{Remark}
\newtheorem{exa}{Example}
\newcommand{\pf}{{\noindent \bf\sl Proof.  }}
\renewcommand{\qed}{\vrule width0pt\hfill \raisebox{-.3ex}
   {\vrule height8pt width8pt depth0pt} \hspace*{-7pt}}
\renewcommand{\Re}{{\bf R}}

\renewcommand{\hbar}{\centerline{\rule{8cm}{0.5mm}}}
\newcommand{\tc}{{\mathcal C}_{\rm tot}}  
\newcommand{\chull}{{\mathcal H}_{\rm cvx}} 
\newcommand{\Acal}{{\mathcal A}}
\newcommand{\Acalh}{{\widehat{\mathcal A}}}
\newcommand{\area}{{\rm Area}}
\newcommand{\ch}{{\widehat C}}  
\newcommand{\kh}{{\widehat k}}

\begin{document}

\large

\title[Soap Film-Like Surfaces Spanning Graphs]  
{\textbf {Area Density and Regularity for Soap Film-Like
 Surfaces Spanning Graphs}}
\author{\bf Robert Gulliver and Sumio Yamada}
\thanks{Supported in part by NSF grant 00-71862}
\date{June 25, 2004}

\begin{abstract}
For a boundary configuration $\Gamma$ consisting of arcs and
vertices, with two or more arcs meeting at each vertex, we treat
the problem of estimating the area density of a
soap film-like surface $\Sigma$ spanning $\Gamma$.  $\Sigma$
is assumed to be strongly stationary for area with respect to
$\Gamma$.  We introduce a notion of total curvature $\tc(\Gamma)$
for such {\em graphs}, or nets,  $\Gamma$.  When the ambient  
manifold $M^n$ has
non-positive sectional curvatures, we show that $2\pi$ times the
area density of $\Sigma$ at any point is less than or equal to
$\tc(\Gamma)$.  For $n=3$, these density estimates imply,
for example, that if $\tc(\Gamma) \leq 22.9 \pi$, then the only
possible singularities of a piecewise smooth  
$(M,\varepsilon,\delta)$-minimizing set $\Sigma$ are curves,  
along which three
smooth sheets of $\Sigma$ meet with equal angles of $120^\circ$.
We also extend these results to allow $M$ to have variable
positive curvature.
\end{abstract}

\maketitle \baselineskip=1.2\normalbaselineskip

\section{Introduction}\label{intro}

The investigation of minimal surfaces has proved extremely fruitful
in a wide range of topics in geometry.  One of the essential
breakthroughs in the subject is the solution of the Plateau
problem by Douglas and by Rad\'o, that is, the construction of a
disc type minimal surface spanned by a Jordan curve $\Gamma$ in
${\bf R}^n$ \cite{D1}, \cite{R}.
Plateau's original motivation was, in part, to study the geometry of soap
films spanned by variously shaped wires. In particular, it is
natural to want to generalize the boundary condition imposed by
Douglas and Rad\'o that the wire $\Gamma$ spanning the surface be a
Jordan curve, or a union of Jordan curves ({\it cf.} \cite{D2},
where $\Sigma$ is a branched immersion of higher topological type).
In this paper, we will
introduce a class of surfaces $\Sigma$ in an ambient manifold $M$,
having a piecewise smooth boundary $\Gamma$ which is    
homeomorphic to a
{\em graph}, that is, a a finite $1$-dimensional polyhedron
(sometimes called a ``net").
Each surface is to
satisfy a regularity condition, and is stationary for
area under variations induced by one-parameter families of
diffeomorphisms of the ambient manifold.  This setting allows us
to consider surfaces whose induced topology is not locally
Euclidean, such as the singular surfaces which may readily
be observed in soap film experiments.  The main
theorems of this paper provide descriptions of the possible
singularities of those minimal surfaces in terms of the geometry
of the boundary set $\Gamma$.

In a Riemannian manifold $M^n$,
we shall consider an embedded graph $\Gamma$ which is a union of
arcs $a_k$ meeting at vertices $q_j$, each of which has
valence at least two.
The {\em valence} of a vertex $q$ is the number of times $q$
occurs as an endpoint of the $1$-simplices $a_k$.
Each $1$-simplex $a_k$ is assumed to be $C^2$,
and to meet its end points with $C^1$ smoothness;  thus there is a
well-defined tangent vector $T_k$ to each $1$-simplex $a_k$ at a
vertex.  At a vertex $q_j$ of valence $d$, we consider the
{\em contribution to total curvature} at $q_j$:
\begin{equation}\label{tceq}
 {\rm tc}(q_j) :=
\sup_{e \in T_{q_j}M}
  \left\{
\sum_{\ell=1}^d \left(\frac{\pi}{2} - \beta_j^\ell(e) \right)
  \right\}
\end{equation}
where $\beta_j^\ell=\beta_j^\ell(e) \in [0,\pi]$ is the   
angle between the tangent vector
$T_\ell$ to $a_\ell$ at $q_j$ and the  
vector $e$.
We define the {\em total curvature} of $\Gamma$ as
\begin{equation}\label{tcdef}
\tc(\Gamma) := \int_{\Gamma^{\rm reg}} |\vec{k}| \,ds +
\sum \{ {\rm tc}(q): q {\rm \ a\ vertex\ of\ } \Gamma\}.
\end{equation}
where $\vec{k}$ is the geodesic curvature vector of $a_i$ as a
curve in $M^n$, and
$\Gamma^{\rm reg} = \Gamma \backslash \{\mbox{vertices} \}$.
It should be noted that our definition of total curvature coincides  
with the standard definition in the case when $\Gamma$ is a
piecewise smooth Jordan curve: the integral of the norm of
geodesic curvature vector plus the sum of the exterior angles
at the vertices.  Namely, in that case, every vertex $q$
of the graph $\Gamma$ is of valence two;  the supremum in         
equation \eqref{tceq} is assumed at vectors $e$ lying in          
the smaller angle between the tangent vectors $T_1$ and $T_2$ to  
$\Gamma$.                                                         
Recall that the {\it density} of $\Sigma$ at $p$ is
$$
\Theta_{\Sigma} (p) := \lim_{\varepsilon \rightarrow 0}
\frac{\Area(\Sigma \cap B_{\varepsilon}(p))}{\pi\varepsilon^2},
$$
provided this limit exists.  The type of surface we will consider in this
paper is a  set $\Sigma$ which is a finite union of $C^2$-smooth
open two-dimensional embedded manifolds $\Sigma_i$, $C^1$ up to
the boundary $\partial \Sigma_i$, with $\partial \Sigma_i$
piecewise $C^1$.  We further impose that the graph $\Gamma$ is a
subset of $S := \cup_i \partial \Sigma_i$.  The class of such surfaces
will be denoted by ${\mathcal S}_{\Gamma}$. Note that given a
surface $\Sigma$ in ${\mathcal S}_{\Gamma}$, the density
$\Theta_{\Sigma}(p)$ is a well defined, upper semi-continuous
function on $\Sigma$.  Moreover, for $\Sigma$ in the class
${\mathcal S}_{\Gamma}$, we may also write
$$
\Theta_{\Sigma} (p) = \lim_{\varepsilon \rightarrow 0}
\frac{{\rm Length}(\Sigma\cap\partial B_{\varepsilon}(p))}{2\pi\varepsilon},
$$
A surface $\Sigma$ in ${\mathcal S}_{\Gamma}$ is said to be {\it
strongly stationary with respect to $\Gamma$} if the first
variation of the area of the surface
is at most equal to the integral over $\Gamma$ of the length of
the component of the variation vector field normal to $\Gamma$
\cite{EWW}.

We can now state the main area-density estimate for the case when
the ambient space is Euclidean
(see Corollary \ref{tc.bound} below): \\

\noindent{{\bf Area-Density Estimate:} {\it Let $\Sigma$  in the
class ${\mathcal S}_{\Gamma}$ be a strongly stationary surface in
$\Re^n$ with respect to its boundary set $\Gamma$. Then
\[
2\pi \Theta_{\Sigma} (p) \leq \tc(\Gamma).
\]
}

\medskip
This estimate is a consequence of two inequalities, the first
being the comparison of area density of $\Sigma$ and of the cone
$C_p(\Gamma)$.  Here, and in the remainder of this paper, for a
point $p \in \Re^n$ and a set $\Gamma_0 \subset \Re^n$, we write
the {\em cone over $\Gamma_0$} as  
$$ C_p(\Gamma_0) := \{p+t(x-p):  x \in \Gamma_0, 0\leq t\leq 1\}.$$
In section 6, where $\Re^n$ is replaced more generally by a
strongly convex Riemannian manifold,  $C_p(\Gamma_0)$ will denote
the geodesic cone over $\Gamma_0$.\\

\noindent{{\bf Theorem~1:}}\,\,{\it Given a strongly stationary
surface $\Sigma$ in ${\mathcal S}_{\Gamma}$, and a point $p$ in
$\Sigma \backslash \Gamma$, let $C_p(\Gamma)$ be the cone spanned
by $\Gamma$ with its vertex at $p$.  Then we have
$$
\Theta_{\Sigma}(p) < \Theta_{C_p(\Gamma)} (p)
$$
unless $\Sigma$ is a cone over $p$ with planar faces.}\\

The second inequality follows from the Gauss-Bonnet formula
applied to the {\it double cover} of the cone $C_p(\Gamma)$.
(We have not found a    
useful Gauss-Bonnet formula for
general $2$-dimensional Riemannian polyhedra in the literature.)\\

\noindent{{\bf Theorem~2} (and {\bf Corollary~5:})}{\it            
$$
2\pi\Theta_{C_p(\Gamma)} (p) =  -\sum_{k=1}^n
\int_{a_k} \vec{k} \cdot \nu_C \,ds + \sum_{k =1}^n \sum_{j =1, 2}
\left(\frac{\pi}{2} - \beta^k_j\right)  \leq \tc(\Gamma),
$$

\noindent                                                         
where $\vec{k}$ is the geodesic curvature vector of $a_k$
in ${\bf R}^n$, $\nu_C$ is the outward unit normal vector at
$a_k \subset \partial C_p (\Gamma)$, and $\beta^k_j$ is the angle
between the tangent vector to $a_k$ at its endpoint $q_j$ and the
line segment from $q_j$ to $p$. }\\                               

The area density estimate
$2\pi\Theta_{\Sigma} (p) \leq \tc(\Gamma)$, when $\Gamma$ is a
rectifiable Jordan curve, is a major ingredient of the work by
Ekholm, White and Wienholtz~\cite{EWW}, where it was proven that
if $\tc(\Gamma) \leq 4 \pi$, then every stationary branched
minimal surface
$\Sigma$ in ${\bf R}^n$ spanned by $\Gamma$ is embedded;  and that
given a compactly supported rectifiable varifold $\Sigma$ which is
strongly stationary with  respect to $\Gamma$ and with area
density $\geq 1$ on $\Sigma \backslash \Gamma$, the
inequality $\tc(\Gamma) < 3 \pi$ implies that $\Sigma$ is smooth
in the interior.  Therefore one can view the results in this paper
as partial extensions of those theorems in~\cite{EWW}, when the
Jordan curve $\Gamma$ of~\cite{EWW} is replaced, more generally,
by a graph.

By imposing appropriate upper bounds on the total curvature   
of the graph $\Gamma$, we obtain the following statements.  We
will denote by $C_Y=3/2$ the area density at its vertex
of the Y-singularity
cone composed of three planes meeting at $120^{\circ}$, and by
$C_T= 6 \cos^{-1} (-1/3) \approx 11.468$ the area density     
at its vertex
of the T-singularity cone spanned by the one-skeleton of the
regular tetrahedron with vertex at its center.  \\

\noindent{{\bf Theorem~3:}}\,\,{\it Suppose $\Gamma$ is a graph in
$\Re^n$ with $\tc(\Gamma) \leq 2\pi C_Y = 3\pi$, and let $\Sigma$
be a strongly stationary surface relative to $\Gamma$ in the class
${\mathcal S}_{\Gamma}$.  Then $\Sigma$ is an embedded surface or
a subset of the Y-singularity cone.} \\
\noindent{{\bf Theorem~4:}}\,\,{\it Suppose $\Gamma$ is a graph in
$\Re^3$ with $\tc(\Gamma) \leq 2 \pi C_T$, and let $\Sigma$ be an
$(M, \varepsilon, \delta)$-minimizing set with respect    
to $\Gamma$ in ${\mathcal S}_{\Gamma}$.  
Then $\Sigma$ is a surface with possibly Y-singularities
but no other singularities, unless it is a subset of
the T-stationary cone, with planar faces.}\\

For the definition od $(M, \varepsilon, \delta)$-minimizing sets,
see Definition \ref{Med} below.  \\[2mm]

In ${\bf R}^3$, there are many known examples of strongly
stationary surfaces.  In particular, there are exactly
ten stationary cones spanned by a graph  
$\Gamma$ on a unit sphere~\cite{AT}.  
Each graph consists of geodesic segments on the sphere meeting in
threes at angles of $120^{\circ}$, including the planar case,
where the graph is simply one great circle spanning ${\bf R}^2$.    
By ordering those ten minimal cones with
respect to the density at the cone vertex, which is the center of
the unit sphere, one has a list of possible tangent cones at the
interior points of an     
$(M, \varepsilon, \delta)$-set $\Sigma$.   The
first three on the list, that is, the ones with the smallest
densities at the vertex, are the plane with its density $1$ where the
graph $\Gamma$ is a great circle;  the Y-singularity cone with its
density $C_Y = 3/2$ where $\Gamma$ consists of  
three semicircles meeting
at the north and south poles at angles of $120^{\circ}$;  and the
T-singularity cone with density $C_T = 6 \cos^{-1} (1/3)$.  
Recall that those ten cones are stationary, but not minimizing,  
under interior deformations.  Hence given one of
those graphs, there may be another surface which is also strongly
stationary  with respect to the same graph, but has strictly
smaller area. Indeed when it comes to soap films, the first three
on the list are the only tangent cones experimentally observed in
the interior of soap films.  This is also true for the
mathematical model in terms of $2$-rectifiable sets, a result
shown by Jean Taylor~\cite{T}:\\[3mm]

\noindent{{\bf Regularity Theorem for Soap Films:}}\,\,            
{\it Away from $\Gamma$, an $(M, \varepsilon, \delta)$-minimizing 
set with respect to $\Gamma$ consists of real
analytic                                  
surfaces meeting smoothly in
threes at $120^{\circ}$ angles along smooth curves, with these curves
in turn meeting in fours at angles of $\cos^{-1}(-1/3)$.}\\[3mm]

The singular curves were proved to be $C^{1, \alpha}$ by
Taylor~\cite{T}, and later shown to be  
real analytic in~\cite{KNS}.  The class
${\mathcal S}_{\Gamma}$ of surfaces we consider in this paper is
chosen so that, given a graph $\Gamma$, we expect to find that
every $(M, 0, \delta)$-minimizing set relative to $\Gamma$
is in the regularity class
${\mathcal S}_{\Gamma}$.  The ten stationary cones
described above are in fact in ${\mathcal S}_{\Gamma}$.
However, due to the lack of understanding of
boundary regularity  of such $(M, 0, \delta)$-minimizing sets, 
it is not yet
known that in general $(M, 0, \delta)$-sets are indeed elements of
the class ${\mathcal S}_{\Gamma}$.

In {\bf Section \ref{non0curv}}, we turn our attention to the case
where the ambient manifold is of variable curvature.  The lack of
homogeneity of the ambient space forces us to consider a
comparison space of constant sectional curvatures, as was done
previously in~\cite{CG2}.  We consider two classes of Riemannian
manifolds $M$ which are strongly convex (not necessarily complete):
manifolds with sectional curvature $K_M$ bounded above by
$- \kappa^2 \leq 0$, and manifolds
with sectional curvature bounded above by $\kappa^2 > 0$.  
For a Euclidean ambient space, as seen above in Theorem 2,    
the area density of the
surface is bounded above by the total curvature of $\Gamma$.
In the variable curvature case, the total curvature of $\Gamma$
is not invariant under diffeomorphisms of $M$ which mimic the
homotheties of $\Re^n$.  Thus, in order to have significance for
both large graphs $\Gamma$ and for small ones, $\tc(\Gamma)$ needs to
be replaced in the following manner:  \\

\noindent {\bf Area-Density Estimate:} ($K_M \leq -\kappa^2$ case)
\,\, {\it Let $\Sigma$ be a strongly stationary surface relative
to $\Gamma$ in the class ${\mathcal S}_{\Gamma}$ in $M_{K \leq
-\kappa^2}$.  Then
$$
2\pi\Theta_p(\Sigma) \leq \tc(\Gamma) - \kappa^2 {\mathcal A}(\Gamma),
$$
where ${\mathcal A}(\Gamma)$ is the minimum cone area of all the
cones with vertex in the convex hull of the set $\Gamma$. } \\

\noindent {\bf Area-Density Estimate:} ($K_M \leq \kappa^2$ case)
\,\, {\it Let $\Sigma$ be a strongly stationary surface relative
to $\Gamma$ in the class ${\mathcal S}_{\Gamma}$ in $M_{K \leq
\kappa^2}$.  Then

$$
2\pi\Theta_p(\Sigma) \leq
\tc(\Gamma) + \kappa^2 \widehat{{\mathcal A}}(\Gamma),
$$
where $\hat{{\mathcal A}}(\Gamma)$ is the maximum spherical
area of all
the cones with vertex in the convex hull of the set $\Gamma$. } \\

We would like to acknowledge fruitful conversations with Brian
White and with Jaigyoung Choe during the development of this
research.

\section{Density and the Regularity of Strongly Stationary Surfaces}
\label{denssec}

Let $\Gamma \subset \Re^n$ be a graph,  
consisting of immersed arcs $a_i$, which are $C^2$ in the interior
and $C^1$ up to their vertices, as in Section \ref{intro};    
we assume that each vertex  has valence at least two.
Let the class ${\mathcal S}_{\Gamma}$ of singular surfaces
be defined as in Section \ref{intro}: for
$\Sigma \in {\mathcal S}_{\Gamma}$, $\Gamma$ is a subset      
of the one-dimensional part $S\subset \Sigma$.                
Within the class ${\mathcal S}_{\Gamma}$, we will look at the
surfaces satisfying the following property.\\

\begin{defn}~\cite{EWW}\label{strongst}  
A rectifiable varifold $\Sigma$ in $\Re^n$ is called {\em strongly
stationary with respect to} $\Gamma$ if for all smooth
$\phi : \Re \times \Re^n$ with $\phi(0,x) \equiv x$, we have
\[
\frac{d}{dt} \Big( \area (\phi(t, \Sigma))+ \area (\phi([0,t] \times
\Gamma)) \Big) \Big|_{t=0} \geq 0.
\]
\end{defn}

The regularity condition on each $\Sigma_i$
guarantees that at almost every point $p$ of $S$, there exists a
unit vector $\nu_{\Sigma_i}$, normal to $\Gamma$, tangential to
$\Sigma_i$, pointing {\it out} of $\Sigma_i$.  Hence on
each $\Sigma_i$ we have the divergence theorem
\[ \int_{\Sigma_i} {\rm div}_{\Sigma_i} X^T dA = \int_{\partial
\Sigma_i} \langle X, \nu_i \rangle ds,  \]
where $X^T$ is the tangential (to $\Sigma_i$) component of $X$.
Note that the strong stationarity condition implies that $\Sigma
\in {\mathcal S}_{\Gamma}$ is stationary, i.e. the first variation
of the area vanishes under deformations supported away from
$\Gamma$:
\[
\frac{d}{dt} \Big( \area (\phi (t, \Sigma)) \Big) \Big|_{t=0} = 0.
\]
If in addition the vector field $X = \phi'(0, x)$ is supported
away from the singular set $S$, then the stationarity condition
implies that the interior of each $\Sigma_i$ is minimal, i.e. its
mean curvature vector $\vec{H}$ vanishes.  

If the vector field $X$ is supported away from $\Gamma$, the
stationarity condition implies
\[
\int_S \sum_{j \in J \subset I} \langle \nu_{\Sigma_j}(p), X^{\perp}(p) \rangle  ds(p) -
\sum \int_{\Sigma_i} \langle \vec{H}, X^{\perp} \rangle dA =0
\]
where $J =J(p)$ indexes the collection of surfaces $\Sigma_j$
which meet at a point $p$ in $S \backslash \Gamma$.  Note that the
second term vanishes since $\vec{H} \equiv 0$.  Since the choice of $X$
is arbitrary, it follows that the vector
\begin{equation}
\nu_\Sigma(p):=\sum_{j \in J(p) \subset I} \nu_{\Sigma_j}(p) = 0
\label{int.balance}
\end{equation}
 almost everywhere on $S \backslash \Gamma$, which we call the
{\it balancing} of $\nu_{\Sigma_i}$ along the singular curves of
$\Sigma$, away from $\Gamma$.  

The strong stationarity condition of a varifold with respect to
$\Gamma$ is equivalent
to the existence
of an ${\mathcal H}^1$-measurable normal (to $\Gamma$) vector field
$\nu$ on $\Gamma$ with $\sup |\nu| \leq 1$ such that
\begin{equation}\label{strstat}
\int_{\Sigma} {\rm div}_{\Sigma} X\, dA =
\int_{\Gamma} \langle X, \nu \rangle ds
\end{equation}
for all smooth vector fields $X$ on $\Re^n$ (see section 7 of
\cite{EWW}).  Note that since $X$ is an ambient vector field along
$\Sigma$,  
${\rm div}_{\Sigma} X$ is the trace on $\Sigma$ of the
ambient covarant derivative of $X$.

In our context, that is, when $\Sigma$ is in
${\mathcal S}_{\Gamma}$, the ${\mathcal H}^1$-measurable vector
field $\nu = \nu_\Sigma$ arises as
\begin{equation}\label{nusig}
\nu_{\Sigma}(p) = \sum_{j \in J(p)} \nu_{\Sigma_j} (p)
\end{equation}
for each $p \in \Gamma$, where $j\in J(p)$ whenever
$p \in \overline{\Sigma_j}$.

\medskip

First we define a surface $\Sigma$ in ${\mathcal S}_{\Gamma}$ to
be {\it locally minimizing relative to $\Gamma$} at $p$ if for a
neighborhood $U$ of $p$, there exists a smaller neighborhood $V$
of $p$ such that for any $\tilde{\Sigma}$, if
$\tilde{\Sigma} \backslash V = \Sigma \backslash V$ and
$\partial \tilde{\Sigma} = \Gamma$, then
${\rm Area} (\tilde{\Sigma}) \geq {\rm Area} (\Sigma)$.  We are
particularly interested in the case $p$ is a point on $\Gamma$.

For intuition, it is useful to understand  the relation between strong
stationarity and the local minimizing property within the class
${\mathcal S}_{\Gamma}$.  First define an ${\mathcal H}^1$
measurable vector field $\nu_{\Sigma}$ defined on $\Gamma$,
as in \eqref{nusig}
above.  The following proposition may be proved using well-known  
methods of the calculus of variations.

\begin{prop}\label{locmin}
Suppose that $\Sigma$ is a surface in ${\mathcal S}_{\Gamma}$.
Then $\Sigma$ is locally minimizing relative to $\Gamma$ at each
point of $\Gamma$ if and only if $|\nu_{\Sigma}| \leq  1$ \,
${\mathcal H}^1$-almost everywhere on $\Gamma$,              
$\nu_{\Sigma} =0$\, ${\mathcal H}^1$-almost everywhere on    
a neighborhood of $\Gamma$ in $S \backslash\Gamma$           
and the regular parts of
$\Sigma$ have vanishing mean curvature vector $\vec{H}$ in some
neighborhood of $\Gamma$.
\end{prop}

\medskip

This proposition says that within the class
${\mathcal S}_{\Gamma}$, the local minimizing property relative
to $\Gamma$ and stationarity away from $\Gamma$ imply
strong stationarity with respect to $\Gamma$.  We remark here that
strong stationarity is strictly weaker than the locally     
area minimizing condition.  In particular, there are surfaces
which are strongly stationary but not locally area minimizing at
certain interior points.  One such example is the cone
$\Sigma \subset \Re^3$ spanned by the $1$-skeleton $\Gamma$ of a
cube, with its vertex at the center of the cube.  It is strongly
stationary relative to $\Gamma$, but is not locally minimizing at
the cone vertex.  Namely, there exists a one parameter family of
polyhedral surfaces of strictly smaller area, in which a
neighborhood of the vertex at the center is replaced by the
2-skeleton of a small cube;  the variation is supported in an
arbitrarily small neighborhood of $p$ \cite{T}.

Next we introduce the following definition, which,
for surfaces in the class ${\mathcal S}_{\Gamma}$, allows us to
isolate the two independent parts of the strong stationarity
condition.  In fact, strong stationarity for surfaces in the class
${\mathcal S}_{\Gamma}$ is equivalent to stationarity in
$\Re^n \backslash \Gamma$ plus the following boundary condition.

\begin{defn}
$\Gamma$ is said to be a {\em variational boundary} of a surface
$\Sigma$ if there exists an ${\mathcal H}^1$ measurable vector
field $\nu_{\Sigma}$ along $\Gamma$ which is orthogonal to
$\Gamma$, with $|\nu_{\Sigma}| \leq 1$ a.e., such
that for all smooth vector fields $X$ defined on $\Re^n$,
$\int_{\Sigma} {\rm div}_{\Sigma} X^{T} \, dA =
\int_{\Gamma} \langle X, \nu_{\Sigma} \rangle \, ds$.
\end{defn}

Observe that Definition \ref{strongst} of strong stationarity  
refers to ambient derivatives of $X$, in contrast with Definition
2, which is intrinsic to $\Sigma$.

\medskip

Now we are ready to state and prove the main result of this section.

\begin{thm}\label{dens.comp}
Given a strongly stationary surface $\Sigma$ in
${\mathcal S}_{\Gamma}$, and a point $p$ of
$\Sigma \backslash \Gamma$, let $C_p(\Gamma)$ be the cone spanned
by $\Gamma$ with its vertex at $p$.  Then we have the following inequality:
\[
\Theta_{\Sigma}(p) < \Theta_{C_p(\Gamma)} (p),
\]
unless $\Sigma$ is a cone over $p$ with planar faces, in which case
we have equality.
\end{thm}

\pf
Let $G(x)$ be the test function    
$\log \rho(x)$, where $\rho(x) = |x - p|$.  $G(x)$ is
the Green's function for the Laplace operator defined on
two-dimensional subspaces of $\Re^n$ which contain the point $p$.
On the other hand, on a minimal surface in $\Re^n$,
the function $G(x)$ is subharmonic, as a
consequence of the trace formula:
\begin{equation}\label{trace}
\triangle_{\Sigma} G =
\sum_{\alpha = 1}^2 \overline{\nabla}^2 G (e_{\alpha}, e_{\alpha}) +
d G (\vec{H}),
\end{equation}
where $\overline{\nabla}$ is the covariant derivative for the
ambient manifold $\Re^n$ (see ~\cite{CG1}).
Thus, we have the following integral estimate:
\[
0 \leq \int_{\Sigma_i \backslash B_{\varepsilon}(p)}
\triangle_{\Sigma_i} G \,\, dA =
\int_{\partial (\Sigma_i \backslash B_{\varepsilon}(p))}
\frac{1}{\rho}\frac{\partial\rho}{\partial\nu_{\Sigma_i}}\,\, ds  
\]
for each $i$, where $\Sigma = \cup_{i\in I} \Sigma_i$ is a surface
in the class ${\mathcal S}_{\Gamma}$.  The equality is due to the
divergence theorem.  Note that each boundary
$\partial (\Sigma_i \backslash B_{\varepsilon}(p)) $
consists of three parts:
\[
\partial (\Sigma_i \backslash B_{\varepsilon}(p)) =
\Big( \partial \Sigma_i \cap \Gamma \Big) \cup
\Big( \partial B_{\varepsilon}(p) \cap \Sigma_i \Big) \cup
\Big( \partial \Sigma_i \cap (S \backslash \Gamma)
 \Big),
\]
since $S = \cup \partial \Sigma_i$.
Now we sum the inequality above over $i$ and reorganize the
boundary terms:
\[
0 \leq \int_{\Gamma} \frac{1}{\rho}
\frac{\partial \rho}{\partial \nu_{\Sigma}} \,\,\, ds
+ \int_{\partial B_{\varepsilon}(p) \cap \Sigma}  \frac{1}{\rho}
\frac{\partial \rho}{\partial \nu_{\Sigma}} \,\,\, ds
+ \sum_{i\in I} \int_{\partial \Sigma_i \cap (S \backslash \Gamma) }
\frac{1}{\rho} \frac{\partial \rho}{\partial \nu_{\Sigma_i}} ds,
\]
where $\nu_{\Sigma}$ is as in equation \eqref{nusig}.  The last term
vanishes, since we have the balancing condition among the unit vectors
$\nu_{\Sigma_i}$ normal to the edges of
$\Sigma_i \cap (S \backslash \Gamma)$,
and tangent to
$\Sigma_i$, pointing outward of $\Sigma_i$, as a consequence of
the (interior) stationarity~(\ref{int.balance}) of $\Sigma$:
%
\[
\sum_{j \in J(p)} \frac{\partial \rho}{\partial \nu_{\Sigma_j}} =
\left\langle\overline{\nabla}\rho,\sum_{j\in J(p)}\nu_{\Sigma_j}\right\rangle
= 0,
\]
for each $p \in S\backslash \Gamma$, where $J(p)$ is the
collection of $j\in I$ with $p\in \overline{\Sigma_j}$.

As for the second term, note that as $\varepsilon$ goes to zero,
$\frac{\partial \rho}{\partial \nu_{\Sigma_i}}$ approaches $-1$
uniformly, and hence
\[
\int_{\partial B_{\varepsilon}(p)\cap\Sigma}\frac{1}{\rho}\frac{\partial\rho}{\partial \nu_{\Sigma}} ds
\]
converges to
\[
\lim_{\varepsilon\rightarrow 0}\left(-\frac{1}{\varepsilon}\right)
{\rm Length}(\Sigma\cap\partial B_{\varepsilon}(p))
= - 2 \pi \Theta_{\Sigma} (p).
\]

Therefore we have obtained the following upper bound for the area
density of $\Sigma$ at $p$:
\begin{equation}\label{sigma.density}
2 \pi \Theta_{\Sigma}(p) \leq
\int_{\Gamma}\frac{1}{\rho}\frac{\partial\rho}{\partial\nu_{\Sigma}} ds.
\end{equation}

We repeat the argument for the surface Laplacian of $G(x)$,
this time replacing $\Sigma$ with  
the cone $C_p(\Gamma)$ spanned by
$\Gamma$ with vertex $p$.                    
Recall that $\Gamma = \cup_j a_j$ where each arc $a_j$ is
$C^2$-regular, $C^1$ up to the end points.  Denote by $A_j$ the
cone $C_p(a_j)$ spanned by $a_j$ with its vertex at $p$.
Thus the cone $C_p (\Gamma)$ is the union of all the {\it fans}
$\overline{A_j} = A_j \cup \partial A_j$.  Observe using
\eqref{trace} that  
away from the vertex $p$, $G(x)$ is harmonic on $A_j$ \cite{CG1}.  
Hence we have
\[
0 = \int_{A_j \backslash B_{\varepsilon} (p)}
\triangle_C G(x) \,\, dA  =
\int_{\partial (A_j \backslash B_{\varepsilon}(p))}
\frac{1}{\rho} \frac{\partial \rho}{\partial \nu_{C}} \,\, ds.
\]
As seen above for $\Sigma$, each boundary
$\partial (A_j \backslash B_{\varepsilon}(p))$
consists of three parts;
we sum the equation above over $j$ and reorganize the boundary
terms, and find:
\[
0=\int_{\Gamma}\frac{1}{\rho}\frac{\partial\rho}{\partial\nu_{C}}\, ds +
\int_{\partial B_{\varepsilon}(p)\cap C}
\frac{1}{\rho}\frac{\partial\rho}{\partial\nu_{C}}\, ds +
\sum_j\int_{C_p(\partial a_j)\backslash B_{\varepsilon}(p)}
\frac{1}{\rho}\frac{\partial\rho}{\partial\nu_{A_j}}\, ds,
\]
where $\nu_C = \nu_{C_p (\Gamma)}$ is defined to be
$\sum_j \nu_{A_j}$, with $\nu_{A_j}$ being the unit vector
normal to the boundary $\partial A_j$, and tangent to the fan
$A_j$, pointing out of $A_j$.

The last term vanishes since the vector $\nu_{A_j}$ and
$\overline{\nabla} \rho$ are perpendicular, which makes
$\partial \rho / \partial \nu_{A_j}$ identically zero on
$C_p(\partial a_j)$.  The second term is equal to \\
$-{\rm Length}(C_p(\Gamma)\cap\partial B_{\varepsilon}(p)) /\varepsilon$
which in turn is equal to $- 2 \pi \Theta_C (p)$,
independent of sufficiently small $\varepsilon >0$.  Therefore we
have obtained
\begin{equation}\label{C.density}
2 \pi \Theta_C(p) =
\int_{\Gamma} \frac{1}{\rho}\frac{\partial\rho}{\partial\nu_{C}}\,\, ds.
\end{equation}

Now observe that $\nu_C$ is the  unit vector normal to $\Gamma$
most closely aligned with the gradient of $\rho$ along $\Gamma$, while
$\nu_{\Sigma}$ is normal to $\Gamma$ with $|\nu_{\Sigma}| \leq 1$,
since $\Gamma$ is a variational boundary of $\Sigma$.  Hence we
have the following inequality:
\begin{equation}\label{rho.deriv}
\frac{\partial \rho}{\partial \nu_{C}} \geq \frac{\partial \rho}{\partial \nu_{\Sigma}}
\end{equation}
almost everywhere along $\Gamma$.  By integrating, we have
\[
\int_{\Gamma} \frac{1}{\rho} \frac{\partial \rho}{\partial \nu_{C}}\,\,\, ds
\geq \int_{\Gamma} \frac{1}{\rho} \frac{\partial \rho}{\partial \nu_{\Sigma}} \,\,\, ds
\]
Combining the inequalities~\eqref{sigma.density},
\eqref{rho.deriv}                     
and the equality~\eqref{C.density},   
we finally get
\begin{equation}\label{densig}
2 \pi \Theta_{\Sigma}(p) \leq \int_{\Gamma} \frac{1}{\rho} \frac{\partial \rho}{\partial \nu_{\Sigma}} \,\,\, ds
\leq \int_{\Gamma} \frac{1}{\rho} \frac{\partial \rho}{\partial \nu_{C}}\,\,\, ds = 2 \pi \Theta_C(p).
\end{equation}

If equality occurs in \eqref{densig}, then
$\triangle_\Sigma G \equiv 0$, and the trace formula
\eqref{trace}, along with a computation of $\overline\nabla^2 G$,
implies that $\overline\nabla\rho$ is tangent to $\Sigma$.  Thus
each two-dimensional face of $\Sigma$ is both a regular minimal
surface and a stationary cone in $\Re^n$, and therefore is part of
a plane passing through $p$.
\qed

\section{Total Curvatures of Graphs}
%
%

Let $\Gamma$ be a graph in $\Re^n$,  
consisting of immersed arcs $a_1, a_2, \dots, a_n$,
which are $C^2$ in the interior and
$C^1$ up to the vertices $q_1, q_2, \dots, q_m$.  
Recall the definition \eqref{tcdef} of total curvature
$\tc(\Gamma)$ of a graph $\Gamma$.
The definition \eqref{tceq} of ${\rm tc}(q)$ for a vertex $q$ of a
graph in a manifold is equivalent to the following for a graph in
$\Re^n$:

\begin{defn} \label{tc}
If $q$ is a vertex of valence $d$ of a graph
$\Gamma \subset \Re^n$,
define the {\em contribution at $q$ to the total curvature of
$\Gamma$} as
$$ {\rm tc}(q) :=
\sup_{p \in \Re^n}
\sum_{\ell=1}^d \left(\frac{\pi}{2} - \beta_\ell(p) \right)
$$
where $\beta_1(p), \dots, \beta_d(p)$ are the interior angles at $q$
which the $d$ edges of $\Gamma$ make with the line segment from $p$.
\end{defn}


The usefulness of these definitions will become clear in section
\ref{GBsec} below;  see esp. Theorem \ref{GB}.

It might be noted that even though the    
geodesic                                  
curvature in $\Sigma$ at a smooth point of $\Gamma$ is given by the
tangential component of                   
the curvature vector of $\Gamma$, there is no such appropriate
vector at a vertex.  This is true already at a vertex of degree
$d = 2$, that is, for a piecewise smooth Jordan curve.

\bigskip

In 
this section, we shall collect some
observations about $\tc(\Gamma)$ for specific cases of a graph
$\Gamma \subset \Re^n$.  These will be used for the examples
below, but will not be referred to in the proofs of the theorems.
As those results are elementary, and some of them previously
known, we include brief proofs for the sake of completeness (see
\cite{MY} and references therein for more general discussion on
minimal network problems.)

\medskip

Consider a vertex $q$ of $\Gamma$ of valence $d$, and let
$T_1, \dots, T_d$ be
the unit tangent vectors to $\Gamma$ at $q$.  For a given point
$p \in \Re^n$, as in Definition \ref{tc},                     
we may write $\beta_\ell(p)$ for the angle between
$T_\ell$ and the line segment from $q$ to $p$.  We shall also (by
abuse of notation:  compare equation \eqref{tceq})            
write this angle as $\beta_\ell(e)$, where $e$
is the unit vector $\frac{p-q}{|p-q|}$.  We write
$e = e_0 \in S^2$ for a point where the sum
$\sum_{\ell=1}^d \left(\frac{\pi}{2} - \beta_\ell(e)\right)$
assumes its maximum value ${\rm tc}(q)$.  Since $e_0$ is also the
minimizer of $\sum_{\ell=1}^d \beta_\ell(e)$, it is the  
spherical {\it Steiner point} of $T_1, \dots, T_d$.  Note that the
existence of $e_0$ follows from compactness of $S^2$.

\subsection{Valence three}

\noindent
\begin{prop}  \label{4pi3}
For all $T_1, T_2$ and $T_3 \in S^2$, there exists
$e \in \{T_1, T_2, T_3\}$ so that
$\beta_1(e) + \beta_2(e) + \beta_3(e) \leq 4\pi/3$.
\end{prop}

\pf $T_1, T_2$ and $T_3$ lie in a small (or great) circle $\gamma$
of $S^2$.  Each spherical distance $d(T_i, T_{i+1})$ ($i = 1,2,3
{\rm \ mod\ } 3$) is less than (or equal to) the length of the smaller
arc of $\gamma$ between $T_i$ and $T_{i+1}$, so their sum is at
most the length of $\gamma$, hence $\leq 2\pi$.  Renumber
$T_1, T_2, T_3$ so that $d(T_2, T_3)$ is the largest of the
three distances, and choose $e = T_1$.  Then $\beta_1(e) = 0$,
while $\beta_2(e), \beta_3(e) \leq \frac{2\pi}{3}$.
\qed

\noindent
\begin{cor}  For any vertex $q$ of valence $d = 3$,
${\rm tc}(q) \geq \pi/6$, with equality if and only if the three
unit tangent vectors $T_1, T_2$ and $T_3$ at $q$ are {\em balanced:}
$T_1 + T_2 + T_3 = 0.$
\end{cor}

\pf
By Proposition \ref{4pi3},
$\sup_e\sum_{\ell=1}^3\left(\frac{\pi}{2}-\beta_\ell(e)\right)\geq
\frac{3\pi}{2}-\inf_i\sum_{\ell=1}^3\beta_\ell(T_i) \geq
\frac{\pi}{6}.$

Now suppose that ${\rm tc}(q) = \pi/6$.  As in the proof of
Proposition \ref{4pi3}, the unit tangent vectors
$T_1, T_2, T_3$ lie on a circle $\gamma \subset S^2$.
But $\beta_2(T_1)+ \beta_3(T_1) =
\sum_{\ell=1}^3\beta_\ell(T_1) \geq
\sum_{\ell=1}^3\beta_\ell(e_0) =
\frac{3\pi}{2} - {\rm tc}(q) = \frac{4\pi}{3}$,
while $d(T_2, T_3) \geq \beta_\ell(T_1)$, $\ell = 2,3$, which
implies that $\gamma$ has length $2\pi$.  Thus $\gamma$ is a great
circle and all of the $d(T_i, T_{i+1}) =\frac{2\pi}{3}$.
\qed

\medskip

In specific situations, it is of interest to compute ${\rm tc}(q)$
exactly, or even to identify the spherical Steiner point $e_0$.
The following lemma is not difficult to prove, using the first
variation of the sum of distances on $S^2$.

\noindent
\begin{lem}\label{3e0}
Suppose a vertex $q$ of $\Gamma$ has valence three, with unit
tangent vectors $T_1, T_2, T_3$ to $\Gamma$ at $q$. Let $e_0$ be a
Steiner point for $T_1, T_2, T_3$.  For $\ell = 1,2,3$ choose a
minimizing geodesic (great circle) in $S^2$ from $e_0$ to   
$T_\ell$, and let $\xi_\ell \in T_{e_0}S^2$ be the unit tangent
vector at $e_0$ to the geodesic.  Then either (1) $\xi_1 + \xi_2 +
\xi_3 = 0$, that is, the geodesics make equal angles $2\pi/3$ at
$e_0$;  or (2) $e_0 = T_\ell$ for some $\ell = 1,2,3$, and the
remaining two vectors $\xi_{\ell + 1}, \xi_{\ell + 2}$ form an
angle $\geq 2\pi/3$ (subscripts {\em modulo} $3$).
\end{lem}

\medskip

For equilateral spherical triangles, one might expect the Steiner
point $e_0$ of the vertices to be the center of the triangle;
however, if the triangle is too large, $e_0$ can only be one of
the corners of the triangle:

\noindent
\begin{cor}\label{equilat}
If the vertex $q$ of $\,\Gamma$ has valence $3$ and its unit tangent
vectors $T_1, T_2, T_3$ make equal angles with each other, then
\begin{equation}\label{tctri}
{\rm tc(q)}=\left\{\begin{array}{ll}
     3\left(\frac{\pi}{2}-\beta \right) & {\rm if} \beta\leq R_0,\\
     \frac{3\pi}{2} - 4\sin^{-1}(\frac12 \sqrt{3} \sin \beta) &
           {\rm if} \beta\geq R_0;
     \end{array}
        \right.
\end{equation}
where $0\leq \beta \leq \pi/2$ is the circumradius, the common
spherical distance from $T_\ell$ to the closer center $N$, of the
triangle formed by $T_1, T_2, T_3$;  and where $R_0 \approx 1.33458$
radians is the value of $\beta$ which makes the two options in
formula \eqref{tctri} equal.
\end{cor}

\pf
It follows from Lemma \ref{3e0}  that a minimizer of
$\sum \beta_\ell$ must be one of the five points
$N, -N, T_1, T_2$ or $T_3$.  But
$\sum \beta_\ell(-N) \geq \sum \beta_\ell(N) = 3\beta$, and
$\sum \beta_\ell(T_i) = 4s$, $i = 1,2,3$, where $2s$ is the side
of the equilateral triangle:  $\sin s = \sin \beta \sin(\pi/3)$.
But $3\beta - 4s$ has the same sign as $\beta - R_0$.

\qed

\subsection{Even valence}

\noindent
\begin{prop}  \label{d=4}
If \, $T_1, T_2, T_3$ and $T_4$ are points on $S^2$, then any of
the Steiner points $e_0$ must be one of the $T_\ell$ or one of the
six (or more) points of intersection of the two great circles
passing through disjoint pairs of the four points $T_\ell$.
\end{prop}

The proof of Proposition \ref{d=4} will be immediate from
the following lemma.  

\noindent
\begin{lem}  \label{opp4}
Let $e_0$ be a Steiner point for $T_1, T_2, T_3, T_4 \in S^2$, and
write $\xi_\ell \in T_{e_0}S^2$ for the initial unit tangent
vector to the minimizing geodesic from $e_0$ to $T_\ell$.  If
$e_0$ is not equal to any of the $T_\ell$, then after reindexing
$\xi_1, \xi_2, \xi_3, \xi_4$ in circular order
around the unit circle of $T_{e_0}S^2$, we have $\xi_1 = -\xi_3$ and
$\xi_2 = -\xi_4$.
\end{lem}

\pf
We compute the first variation of $\sum_{\ell=1}^4 \beta_\ell(e)$,
and find that
$0=-\sum_{\ell=1}^4 \langle \xi_\ell,\xi \rangle$ for any
$\xi\in T_{e_0}S^2$.                             
 We conclude that the $\xi_\ell$ are balanced:
\begin{equation}\label{bal}
\xi_1 + \xi_2 + \xi_3 + \xi_4 = 0.
\end{equation}
Write $\eta_\ell$ for the oriented angle from $\xi_\ell$ to
$\xi_{\ell+1}$, $\ell$ {\em modulo} $4$, with
$0 \leq \eta_\ell \leq 2\pi$.

If $\xi_1 = -\xi_3$, then also $\xi_2 = -\xi_4$ according to
\eqref{bal}, and we are done.  Otherwise, the sum
$\xi_1 + \xi_3$ makes the oriented angle
$\frac12 (\eta_1+\eta_2)$ {\em modulo} $\pi$ with $\xi_1$, while
the sum $\xi_2 + \xi_4$ makes the angle
$\frac12 (\eta_2+\eta_3)$ {\em modulo} $\pi$ with $\xi_2$.  But
$\xi_2 + \xi_4 = -(\xi_1 + \xi_3)$, hence
$\eta_1 + \frac12 (\eta_2+\eta_3) = \frac12 (\eta_1+\eta_2) + \pi$
{\em modulo} $\pi$,                            
implying that $\eta_1 + \eta_3 = 0$ {\em modulo} $2\pi$.  
But $\eta_1 + \eta_2 + \eta_3 + \eta_4 = 2\pi$ and
$\eta_\ell \geq 0$, so this forces
either $\eta_1 = \eta_3 = 0,$ implying         
$\xi_1 = \xi_2$ and $\xi_3 = \xi_4$;  or       
$\eta_2 = \eta_4 = 0,$ implying $\xi_2 = \xi_3$ and $\xi_4 = \xi_1$.  
The conclusion now follows from equation \eqref{bal}
in this case as well.
\qed

\medskip

The following lemma has a complex statement but a straightforward
demonstration.

\noindent
\begin{lem}  \label{sum}
Let $\widetilde{\Gamma}$ and $\widehat{\Gamma}$ be graphs with a
common vertex $\widetilde{q} = \widehat{q}$.  Write ${\Gamma}$ for
the union of $\,\widetilde{\Gamma}$ and $\,\widehat{\Gamma}$, and
write $q\,$ for the common vertex when considered as a vertex of
$\,\Gamma$.  Write $\{\widetilde{T}_1, \dots, \widetilde{T}_k\}$ for the
unit tangent vectors to $\widetilde{\Gamma}$ at $\widetilde{q}$,
and let $\{\widehat{T}_1, \dots, \widehat{T}_{d-k}\}$ be the
unit tangent vectors to $\widehat{\Gamma}$ at $\widehat{q}$.  
Then
${\rm tc}(q) \leq {\rm tc}(\widetilde{q}) + {\rm tc}(\widehat{q})$.
If further
$\{\widetilde{T}_1, \dots, \widetilde{T}_k\}$ and
$\{\widehat{T}_1, \dots, \widehat{T}_{d-k}\}$ share the same Steiner point
$\widetilde{e}_0 = \widehat{e}_0$, then the Steiner point $e_0$ of
$\{\widetilde{T}_1, \dots, \widetilde{T}_k, \widehat{T}_1, \dots,
\widehat{T}_{d-k}\}$ is equal to both, and
${\rm tc}(q) = {\rm tc}(\widetilde{q}) + {\rm tc}(\widehat{q})$.
\end{lem}

\noindent
\begin{cor}  \label{opp}
If a vertex $q$ of $\Gamma$ has an even valence $d$ and the tangent
vectors at $q$ occur in antipodal pairs, then ${\rm tc}(q) = 0$.
\end{cor}

\pf
Observe that a vertex $q$ of degree $2$ in a straight edge, that
is, with $T_2 = -T_1$, has ${\rm tc}(q) = 0$, with {\em any} point
of $S^2$ as a Steiner point.  The conclusion then follows from
Lemma \ref{sum} by induction on $d/2$.
\qed

\medskip

In contrast with Corollary \ref{equilat}, even valence makes
computations easier:

\noindent
\begin{cor}\label{Neven}
For a regular polygon in $S^2$ with an {\em even} number $d$ of
sides, the closer center in $S^2$ of the polygon is a Steiner
point of the corners $T_1, \dots, T_d$.
\end{cor}

\pf
Let $N$ be the closer center (closer than $-N$) of the regular
polygon of $d =: 2k$ sides, with vertices $T_1, \dots, T_{2k}$ in
order.  A Steiner point of two opposite vertices
$\{T_i, T_{k+i}\}$ is any point along the minimizing geodesic arc
joining them, in particular the midpoint $N$.  Now apply Lemma
\ref{sum} via induction on $k$.
\qed

\noindent
\begin{prop}\label{iffeven}
For a vertex $q$ of a graph $\Gamma \subset \Re^3$ with unit
tangent vectors $T_1, \dots, T_d$ all lying in a plane through $0$
and making equal angles, an orthogonal unit vector $N$ is a
Steiner point if and only if $d$ is even.
\end{prop}

\pf
If $d$ is even, the conclusion is given by Corollary
\ref{Neven}.  If $d = 2k+1$ is odd, then the sum
$\sum_{\ell=1}^d \beta_\ell(e)$ equals $(2k+1)\pi/2$ for $e = N$,
and equals $(1+2+\dots+k)4\pi/(2k+1)$ for $e = T_1$, which is
smaller by a difference of $\frac{\pi}{2(2k+1)}$.  Thus $N$
cannot be the Steiner point.  
\qed

\section{Gauss-Bonnet formula for Cones}\label{GBsec}
In this section we will prove a Gauss-Bonnet formula for
two dimensional cones in ${\bf R}^n$.  First we quote the following
classical result.\\

\noindent {\bf Euler's Theorem} (see\cite{O})
{\it For a connected graph $\Gamma'$ with even valence at each  
vertex, there is a continuous mapping of the circle to $\Gamma'$
which traverses each edge exactly once.}\\

An immediate consequence of this result is that {\it any} connected
finite graph $\Gamma$ has a continuous mapping of the circle which
traverses each edge exactly {\it twice.}  Namely, we may apply Euler's
theorem to the graph $\Gamma'$ obtained from $\Gamma$ by doubling
each edge and leaving the vertices alone.  Note that the new graph
$\Gamma'$ has even valence at each vertex.

We shall derive the density formula of Theorem \ref{GB} below in three
steps, beginning from a well known case.

Suppose first that $\Gamma_0$ is a {\bf smooth closed curve} in $\Re^n$,
not necessarily simple, and $p$ a point not on $\Gamma_0$. Without
loss of generality (after a suitable scaling centered at $p$), we
may assume that $\Gamma_0$ lies outside the unit ball $B_1(p)$
centered at $p$.

Define $\Pi_p$ to be the radial projection to the unit sphere
centered at $p$:
\[
\Pi_p : \Re^n \backslash \{p\} \rightarrow \partial B_1(p);
\]
\[
\Pi_p(x) = p + \frac{x-p}{|x-p|}.
\]

Let $A  = C_p(\Gamma_0) \backslash B_1 (p)$ be the annular region
between $\Gamma_0$ and $\Pi_p \Gamma_0$. By the Gauss-Bonnet
formula, we have
\begin{equation}\label{GBform}
-\int_{\partial A} \vec{k} \cdot \nu_C \,ds + \int_A K\, dA =
 2 \pi \chi (A)
\end{equation}
where $\vec{k}$ is the curvature vector of the graph $\partial A$
in $\Re^n$, $\nu_C$ is the outward normal to $\partial A$, $K$ is the
Gauss curvature of $A$, and $\chi (A)$ is the Euler
characteristic of $A$.  For $A$, $K \equiv 0$ and $\chi (A) = 0$.
Hence
\begin{eqnarray*}
0 & = & \int_{\partial A} \vec{k} \cdot \nu_C \,ds \\
  & = & \int_{\Pi_p \Gamma_0} \vec{k} \cdot \nu_C \,ds +
    \int_{\Gamma_0} \vec{k} \cdot \nu_C \,ds
\end{eqnarray*}
For $q \in \Pi_p \Gamma_0$, $\vec{k}(q)$ is          
the unit vector from $q$ to $p$, so that             
the first integral on the last line is equal to  
the length of $\Pi_p \Gamma_0$,
which is also equal to
$2 \pi \Theta_{C_p(\Gamma_0)} (p)$.  Therefore we have for the cone
$C_p({\Gamma_0})$ the following equation:
\begin{equation}
2 \pi \Theta_{C_p(\Gamma_0)} (p) =
{\rm Length} (\Pi_p \Gamma_0) =
-\int_{\Gamma_0} \vec{k} \cdot \nu_C \,ds, \label{GB.smooth}  
\end{equation}
where $\nu_C(q)$ is the unit normal vector to $\Gamma_0$ in the
plane spanned by the tangent vector at $q$ and the vector $p-q$,
and pointing  away from the cone vertex $p$.  Note that
$C_p \backslash \{ p \}$ is flat with respect to the induced
metric, that is, locally isometric to $\Re^2$.  Note further
that the integrand $\vec{k} \cdot \nu_C$ is the {\it intrinsic}
geodesic curvature of $\Gamma_0$ considered as a locally embedded
curve in $C_p$.

Next, when $\Gamma'$ is a {\bf piecewise smooth}
immersion of the circle,  
we generalize the  
formula above as follows.  Let $\Gamma'$ be a union of smooth
segments $a_i$, each of which is $C^2$ in the interior
and $C^1$ up to the end points $q_{i,0}, q_{i,1}$. We denote
$q_{i, j} \sim q_{i', j'}$    if they represent the same point
where $a_i$ and $a_{i'}$ meet.  Then the cone $C_p(\Gamma')$ can
be thought as a union of {\it fans} $A_i(p) = C_p(a_i)$,  
which is the part of the cone $C_p (\Gamma')$ spanned by $a_i$, with
radial edges $\overline{pq_{i,0}}$ and $\overline{pq_{i,1}}$.  The
right hand side of the equation~(\ref{GB.smooth}) then generalizes
as
\begin{equation}
2 \pi \Theta_{C_p(\Gamma')} (p) =
{\rm Length} (\Pi_p \Gamma') =
-\sum_i \int_{a_i} \vec{k} \cdot \nu_C \,ds  
+ \sum_i \sum_{j=1, 2} \left(\frac{\pi}{2} - \beta^i_j\right)
\label{GB.general}
\end{equation}
where $\beta^i_j$ is the angle between $a_i$ and
$\overline{pq_{i, j}}$ as they meet at $q_{i, j}$.  To see how the
last term arises, suppose now that $a_i$ and $a_{k}$  
are the consecutive edges in $\Gamma'$ joined at $q_{i,j} \sim
q_{k, j'}$.  Then the quantity
$(\pi/2 - \beta^i_j) + (\pi/2 - \beta^{k}_{j'}) =
\pi -(\beta^i_j + \beta^{k}_{j'})$ is the amount    the curve
$a_i \cup a_{k}$ turns at $q_{i,j} \sim q_{k,j'}$, when
considered as a locally isometrically embedded curve in $\Re^2$.

Finally, coming back to the {\bf original graph} $\Gamma$, Euler's
theorem says that the graph $\Gamma$ with each edge traced twice
while its vertices are left intact, which we denoted by $\Gamma'$,
can be parameterized by a copy of $S^1$.  Write $\Gamma'$ as the union
of $a'_i$ where each $a_k \,\, (k=1, \dots n)$  arises twice as
$a'_i \,\, (i=1, \dots, 2n)$, as one goes {\it around} $\Gamma'$
once.

Applying the generalized equation~(\ref{GB.general}) when $\Gamma'$
is $\cup_i^{2n} a'_i$, we obtain the following
description of the density of the cone $C_p (\Gamma)$ at $p$.

\noindent
\begin{thm} \label{GB}
With the notations as above we have the following,
\begin{equation}\label{GBeq}
2 \pi \Theta_{C_p (\Gamma)} (p) =
-\sum_{k=1}^n \int_{a_k} \vec{k} \cdot \nu_C \,ds
+\sum_{k =1}^n\sum_{j=1,2}\left(\frac{\pi}{2}-\beta^k_j\right).
\end{equation}
\end{thm}

\pf
>From the preceding discussion, we have
\begin{equation}
2 \pi \Theta_{C_p(\Gamma')} (p) =
{\rm Length} (\Pi_p \Gamma') =
-\sum_{i=1}^{2n} \int_{a'_i} \vec{k} \cdot \nu_C \,ds  
+ \sum_{i=1}^{2n}\sum_{j=1,2}\left(\frac{\pi}{2}-\beta'^i_j\right).
\label{GB.gamma'}
\end{equation}

Note that the length of $\Gamma'$ is twice the length of $\Gamma$.
Also note that when the edges $a'_{i_1}$ and $a'_{i_2}$ of
$\Gamma'$ represent the same edge $a_k$ of $\Gamma$, we have
\[
\int_{a_k} \vec{k} \cdot \nu_C \,ds =           
\int_{a'_{i_1}} \vec{k} \cdot \nu_C \,ds =
\int_{a'_{i_2}} \vec{k} \cdot \nu_C \,ds
\]
independent of the orientations imposed by the Euler circuit.
Lastly, over the whole circuit $\Gamma'$, the quantity $\pi/2 -
\beta_i^j, \,\, (i=1, \dots, n; \,\, j =1,2)$ appears twice.  The
statement of the theorem then follows by dividing both sides of the
equation~\eqref{GB.gamma'} by two.
\qed

\section{Regularity of Stationary Surfaces}

Using the notations from section \ref{denssec} above, we have the
following immediate consequence to {\bf (1)} the density comparison
(Theorem~\ref{dens.comp}) between the area density of a strongly
stationary surface
$\Sigma$ with respect to $\Gamma$ and that of the cone $C_p(\Gamma)$ over
$\Gamma$ with vertex $p$;  and {\bf (2)} the Gauss-Bonnet
formula~(Theorem~\ref{GB}), which estimates the density of the cone in
terms of the total curvature of the graph $\Gamma$:

\noindent
\begin{cor} \label{tc.bound}
The following inequality holds between the area density of a
strongly stationary surface $\Sigma$
and the total curvature $\tc$ of $\Gamma$:
\[
2 \pi \Theta_{\Sigma}(p) \leq \tc(\Gamma).
\]
\end{cor}

\pf
We need only observe that in the conclusion of Theorem \ref{GB},
the right-hand side of equation \eqref{GBeq} is bounded above by
$\tc(\Gamma)$.
\qed

\noindent
\begin{thm}\label{noY}
Suppose $\Gamma$ is a graph in $\Re^n$ with
$\tc(\Gamma) \leq 2\pi C_Y = 3\pi$,
and let $\Sigma$ be a strongly stationary surface relative to
$\Gamma$ in the class ${\mathcal S}_{\Gamma}$.  Then $\Sigma$ is an
embedded surface or a subset of the Y singular cone.
\end{thm}
\pf
At a point $p$ on $\Sigma$, the proof ofthe  above
Corollary~\ref{tc.bound} to the Gauss-Bonnet formula says that
\[
\Theta_{\Sigma}(p) \leq \Theta_{C_p (\Gamma)} (p)
\leq \frac{1}{2 \pi} \tc (\Gamma) \leq C_Y,
\]
where the last inequality is the hypothesis.  If
$\Theta_{\Sigma}(p) < C_Y,$ we claim that
$\Sigma$ is regular at $p$ by the proof of Theorem 7.1
of~\cite{EWW}.  For the sake of completeness, we reproduce their
argument here.

Let $T_p \Sigma$ be the tangent cone at $p$, whose existence and
uniqueness is guaranteed by the regularity assumption we impose on
the class of surfaces ${\mathcal S}_{\Gamma}$. Then
$\Theta_{T_p \Sigma} (x) < 3/2$ for all $x$ in the cone since in
any minimal cone, the highest density occurs at the vertex.  This
is because the density function $\Theta_{T_p \Sigma} (x)$ is upper
semi-continuous (\cite{Si} \S17.8) and constant along radial
lines.  Now the intersection of $T_p \Sigma$ with the unit sphere
is a collection of geodesic arcs~\cite{AA}, which means that the
cone is a polyhedron.  At most two faces of the polyhedron
$T_p \Sigma$ can meet along a radial edge, since otherwise the
density at points along the edge would be $\geq 3/2$.  This means
$T_p \Sigma \cap S^{n-1}$ is a union of complete great circles.
Since the density is $< 3/2$, there is only one great circle and
it has multiplicity $1$.  By Allard's regularity theorem
(\cite{Al} or~\cite{Si}), this means that $\Sigma$ is regular at
$p$.

On the other hand, if $\Theta_{\Sigma}(p) = C_Y,$ then equality
holds in Theorem \ref{dens.comp}, implying that $\Sigma$
itself is a cone with vertex $p$ and planar faces.  But the Y cone
is the unique (up to rotation in $\Re^n$) stationary cone having
density $3/2$.
\qed

\medskip

As seen above, $3/2$ is the first nontrivial upper bound for the
area density above $1$, for the class of surfaces we are studying.
As for a larger upper bound, we will restrict our attention to the
case when the ambient Euclidean space is $\Re^3$.  There are
exactly ten stationary  
cones in $\Re^3$ ~\cite{AT}, where a cone is stationary  
when its intersection with the unit sphere
is a net of geodesics meeting in threes at $120^\circ$.
Ordered with respect to the area density $\Theta$ at the vertices
of the cones, the first three on the list are the plane with
$\Theta =1$;  Y = three half-planes meeting at $120^{\circ}$ with
$\Theta = C_Y = 3/2$; and the cone T spanned by the regular
tetrahedron with
$\Theta = C_T = 6 \cos^{-1}(-1/3) \approx 11.4638$.

In order to state the next result, we need to introduce the following
definition~\cite{Alm}.

\begin{defn}\label{Med}
Let $\varepsilon$ be a bound of the form  
$\varepsilon(r) = Cr^\alpha$ for some $\alpha > 0$, and choose  
$\delta >0$.  We define $\Sigma \subset \Re^n$ to be an  
$(M,\varepsilon,\delta)$-{\em minimal set}  
with respect to $\Gamma \subset \Re^n$  
if $\Sigma$ is $2$-rectifiable and if, for every Lipschitz  
mapping $\Phi: \Re^n\rightarrow \Re^n$ with the diameter $r$ of
the support $W$ of $\Phi -${\rm id} less than $\delta$,  
$$ \mathcal{H}^2(S\cap W) \leq  
\left(1+\varepsilon(r)\right)\mathcal{H}^2\Big(\Phi(S\cap W)\Big).$$  
\end{defn}

We have the following regularity statement in $\Re^3$ for $\Gamma$
with small total curvature.

\noindent
\begin{thm}\label{noT}
Suppose $\Gamma$ is a graph in $\Re^3$ with
$\tc(\Gamma) \leq 2 \pi C_T$, and let $\Sigma$ be an
$(M, 0, \delta)$-minimal surface with $\Gamma$ as its variational
boundary in ${\mathcal S}_{\Gamma}$.  Then $\Sigma$ is a surface
with possibly Y singularities but no other singularities, unless
it is a subset of the T stationary cone, with planar faces.
\end{thm}

\pf
As in the proof of the previous theorem, for each point $p$ in
$\Sigma$, we have a series of inequalities
\[
\Theta_{\Sigma}(p) < \Theta_{C_p (\Gamma)} (p) \leq
\frac{1}{2\pi} \tc (\Gamma) \leq C_T,
\]
unless $\Sigma$ is a cone over $p$ with planar faces.  We now use
results in~\cite{T}(II.2 and II.3), which imply that the tangent
cone of an $(M, 0, \delta)$-minimal set $S$ at $p$ is area-minimizing
with respect to the intersection with the unit sphere centered at
$p$, and that the plane, the Y-cone and the T-cone are the only
possibilities for the tangent cone.  The inequality above implies
that the tangent cone $T_p \Sigma$ can only be the plane or the Y
singularity, since all other stationary singular cones have higher
density.  If there is a point $p$ where the tangent cone to
$\Sigma$ is any other cone than the plane or Y, then it can only
be the T stationary cone.  But in this case,
$\Theta_{\Sigma}(p) = C_T$, and $\Sigma$ itself is a cone over
$p$.  It follows that $\Sigma =$ T.
\qed

\begin{ram}
A surface $\Sigma$ in the class ${\mathcal S}_{\Gamma}$ which is
$(M, 0, \delta)$-minimal with $\Gamma$ as its variational boundary
is in particular strongly stationary with respect to $\Gamma$
(See the remark preceding {\bf Definition 2}.) However note that
a cone over the one-skeleton $\Gamma$ of the cube is strongly
stationary w.r.t. $\Gamma$, but is not an $(M, 0, \delta)$-minimal
set.
\end{ram}

\begin{ram}
The previous papers~\cite{EWW} and \cite{CG2} had consequences for
the knot class of a curve in a 3-dimensional manifold satisfying
an inequality on its total curvature.  Similar consequences for
the isotropy class of a graph would follow from Theorems~\ref{noY}
and~\ref{noT} if the boundary regularity of an area-minimizing
rectifiable set bounded by a graph could be proved.
\end{ram}

\medskip

\noindent
\begin{exa}\label{football}
In this example, we shall show that the hypothesis $\tc(\Gamma)
\leq 3\pi$ of Theorem \ref{noY} is sharp.  Specifically, we shall
construct a graph $\Gamma$ in $\Re^3$ with $\tc(\Gamma) = 3\pi$,
such that a subset of the minimal cone $Y$, including a nonempty
segment of the singular line, is strongly stationary with respect
to $\Gamma$.
\end{exa}

Recall the description of $Y$ in Section 2 above:  $Y$ consists of
three half-planes $P_1, P_2, P_3$ meeting along a line $S$, and
making equal angles $2\pi/3$ at each point of $S$.  Recall also
the angle $R_0 = 1.33458$ radians $= 76.466^o$ of Corollary
\ref{equilat}.

We choose two points $q^\pm$ along $S$, and construct $\Gamma$ as
the union of three $C^2$ convex plane arcs $a_\ell$, where $a_\ell$
joins $q^-$ to $q^+$ in the half-plane $P_\ell$, $\ell = 1,2,3$,
all making an angle
$\alpha^\pm$ with $S$ at the endpoint $q^\pm$, where
$0 < \alpha^\pm \leq R_0$.  Since $a_\ell$ is a convex plane arc,
the integral of $|\vec{k}|$ along $a_\ell$ equals
$\alpha^+ + \alpha^-$.
Using Corollary \ref{equilat}, we
may compute that the contribution at $q^\pm$ to the total
curvature of $\Gamma$ is ${\rm tc}(q^\pm) = 3(\pi/2 - \alpha^\pm)$.
Thus $\tc(\Gamma) =
3(\alpha^+ + \alpha^-) + 3(\pi/2 - \alpha^+) +
3(\pi/2 - \alpha^-) = 3\pi$, as claimed.
\qed

\medskip

In Example \ref{football},  
intuition might lead the reader to expect that every case, with a
skinny or fat angle, would give rise to a sharp inequality.  In
fact, for the case $\alpha^\pm > R_0$,  
the inequality is {\em not} sharp, as follows using Corollary
\ref{equilat}.

\medskip

\noindent
\begin{exa}\label{tetrahedron}
In this example, we shall show that the hypothesis
$\tc(\Gamma) \leq 2\pi C_T$ of Theorem \ref{noT} is sharp.  In
fact, we shall show that the cone $\Sigma$ over     
the one-skeleton $\Gamma$ of the regular
tetrahedron itself provides an example.
\end{exa}

Let $\alpha_T$ be the angle between an edge $\overline{q_k q_i}$  
of $\Gamma$ and $\overline{q_k p}$, $1\leq k<i \leq 4$, where     
$p$ is the center of the tetrahedron.                             
Then $\cos(\alpha_T) =
\sqrt{2/3}$, so $\alpha_T = 0.61548$ radians, which is less than
$R_0 = 1.33458$ radians.  This shows, using Corollary
\ref{equilat}, that $\tc(\Gamma) = 6\pi - 12 \alpha_T$.

On the other hand, we may apply Theorem \ref{GB} above to compare
the total curvature of $\Gamma$ with the density of $\Sigma$ at
the interior singular point $p$.  
Namely, by
Corollary \ref{equilat}, $p$ will be a Steiner point for the unit
tangent vectors at each of the four vertices, and the curvature
vector $\vec{k} \equiv 0$  along the regular part of $\Gamma$.
In the notation of Theorem \ref{GB}, all twelve of the interior
angles $\beta_k^j$, $1 \leq k \leq 6, \ j = 1,2$ are equal to $\alpha_T$.  
Therefore the density $2\pi C_T$ of the cone at $p$ equals
$$\sum_{k =1}^6 \sum_{j =1, 2} \left(\frac{\pi}{2} - \alpha_T\right)
= 6\pi - 12\alpha_T = \tc(\Gamma).$$
\qed

\medskip



Example~\ref{football} illustrates that the upper bound $3 \pi$
for $\tc(\Gamma)$ is achieved for a non-Jordan curve $\Gamma$.
The next proposition in turn says that
among all the embedded graphs $\Gamma$ which are homeomorphic
to the graph of Example~\ref{football}, $3 \pi$ is
the sharp lower bound for the total curvature $\tc(\Gamma)$.

\begin{prop}  \label{3pi}
Let $\Gamma$ be an embedding into $\Re^3$ of the topological graph
with exactly two vertices $q^\pm$ and three edges $a_1,\, a_2$
and $a_3$, each of which has endpoints $q^+$ and $q^-$.  
Then
$\tc(\Gamma) \geq 3\pi$.  Moreover, equality holds if and only if
each $a_\ell$ is a convex plane arc with unit tangent vectors
$T_\ell^\pm$ at $q^\pm$ satisfying the condition that
$\pm e := \pm \frac{q^- - q^+}{|q^- - q^+|}$ is a
Steiner point for the three points $T_1^\pm, T_2^\pm, T_3^\pm$ on
$S^2$, at both $q^-$ and $q^+$.
\end{prop}

\pf
The ``if" part of the equality conclusion follows essentially from
the discussion of Example~\ref{football} above.  We have adapted
the notation introduced there;  further, let $\alpha_\ell^\pm$
be the angle between $T_\ell^\pm$ and the unit tangent vector
$\pm e$ at $q^\pm$ to the closed line segment $L$ joining $q^\pm$
to $q^\mp$.  Then $a_\ell \cup L$ is a closed curve in $\Re^3$, so
by Fenchel's theorem  
$$2\pi \leq \tc(a_\ell \cup L) = \int_{a_\ell} |\vec{k}| \, ds +
(\pi - \alpha_\ell^+) + (\pi - \alpha_\ell^+).$$
Thus
$\int_{a_\ell} |\vec{k}| \, ds \geq \alpha_\ell^+ + \alpha_\ell^-,$
with equality if and only if $a_\ell$ is a convex planar arc.

Meanwhile, ${\rm tc}(q^\pm) :=
\sup_p \sum_{\ell=1}^3 \left(\frac{\pi}{2} - \beta_\ell^\pm(p)\right) \geq
\sum_{\ell=1}^3 \left(\frac{\pi}{2} - \alpha_\ell^\pm\right).$
Further, equality holds if and only if $\pm e$ is a Steiner point
on $S^2$ for the three points $T_1^\pm, T_2^\pm, T_3^\pm$.
Therefore,
\begin{eqnarray*}
\tc(\Gamma) &:=& \sum_{\ell=1}^3 \int_{a_\ell} |\vec{k}| \, ds +
{\rm tc}(q^+) + {\rm tc}(q^-) \\
&\geq& \sum_{\ell=1}^3 \left[(\alpha_\ell^+ + \alpha_\ell^-)+
(\frac{\pi}{2} - \alpha_\ell^+)+ (\frac{\pi}{2} - \alpha_\ell^-)\right] =
3\pi,
\end{eqnarray*}
with equality if and only if $a_\ell$ is a convex planar
arc and $\pm e$ is the Steiner point.
\qed \\

There is a second combinatorial structure for a connected graph $\Gamma$
with two trivalent vertices and three edges:  the ``handcuff"
consisting of two loops plus an arc joining the vertices of the
loops.  Similarly to Proposition \ref{3pi}, it may be shown that
an embedding of such $\Gamma$ in $\Re^3$ must have total curvature
at least $3\pi$.  In fact, it appears likely that the hypothesis
of Theorem \ref{noY} can hold strictly only for the embedded circle
or the two-leafed rose, that is, two circles connected at a point.\\

The next example will be much more complex than those above.

\noindent
\begin{exa}\label{11net}
In this example, we shall construct a graph $\Gamma$ with
$\tc(\Gamma) = 44\pi < 2\pi C_T$, which is sufficiently
complicated that the presence of a $T$-singularity in a strongly
stationary surface $\Sigma$ might appear likely without Theorem
\ref{noT} above.
\end{exa}

Let $\Gamma$ be the union of eleven congruent (convex) plane
ovals. 
$\Gamma$ will consist of six
horizontal copies in planes $\{z = c_k\}$, $1 \leq k \leq 6$,
obtained from each other by translation in the $z$-direction; and
five copies in vertical planes $\{y = c_k\}$, $7 \leq k \leq 11$,
obtained from each other by translation in the $y$-direction.  We
also assume that each vertical oval meets each horizontal oval
twice. 
For clarity, we assume
that each of the eleven ovals includes two unit line segments
tangent to the faces $\{x=0\}$ and $\{x=1\}$ of the unit cube.
In particular, we assume $0<c_1 < c_2< \dots <c_6<1$ and
$0<c_7 < c_8< \dots <c_{11}<1$.

        Then $\Gamma$ has $60$ vertices $q_1, \dots, q_{60}$,
each of valence $d=4$, and
at each vertex, the unit tangent vectors $T_1, T_2, T_3, T_4$
satisfy $T_3 = -T_1$ and $T_4 = -T_2$.  It follows from
Corollary \ref{opp} that ${\rm tc}(q_i) = 0$, $1\leq i \leq 60$.
Each of the eleven ovals contributes $2\pi$ to the total curvature
of $\Gamma_{\rm reg}$.  Therefore $\tc(\Gamma) = 44\pi < 2\pi C_T$.
\qed

\medskip

\section{Nonzero ambient curvature}\label{non0curv}

In this section, we shall indicate the modifications which need to
be made to generalize Theorems \ref{dens.comp}, \ref{GB},
\ref{noY} and \ref{noT} above to the case where the ambient
space $\Re^n$ is replaced by a manifold $M^n$ having  
variable sectional curvatures.  
In the case of an immersed minimal surface (or a branched
immersion) with smooth boundary, the proof was carried out in
\cite{CG2}; the conclusions in subsection \ref{pos}, however, are
more general than those of \cite{CG2}, even in the case of a
Jordan curve $\Gamma$, since \cite{CG2} requires constant curvature
in the positive case.  Many, although not all, of the proofs of
\cite{CG2} can be adapted with little change to the present
context of singular minimal surfaces which are strongly stationary
with respect to a graph $\Gamma$.

For the rest of this section, let $M^n$ be a strongly convex
Riemannian manifold having sectional curvatures bounded above by
either (\ref{neg}) a non-positive constant $-\kappa^2$;  or
(\ref{pos}) a positive constant $\kappa^2$.  $M^n$ is said to be
{\em strongly convex} if any two points are connected by a unique
minimizing geodesic.  For example, $M^n$ might be a complete,
simply connected Hadamard-Cartan manifold, or a convex open subset
of such a complete manifold, or a convex open subset of
a ball of radius $\pi/\kappa$ in a complete, simply connected
manifold $M^n$ with sectional curvatures $K_M\leq\kappa^2$.

\subsection{Nonpositively Curved Manifold}\label{neg}

Let $M^n$ be a strongly convex Riemannian manifold whose sectional
curvatures are bounded above by a non-positive constant
$-\kappa^2$.  We consider a graph $\Gamma \subset M^n$ and a
surface $\Sigma$ in the class ${\mathcal S}_{\Gamma}$ which is
strongly stationary with respect to $\Gamma$.

Choose a point $p$ of $\Sigma$.  We shall {\em assume} that
$\Gamma$ is nowhere tangent to the minimizing geodesic from $p$;
the general cases of Theorems \ref{dens.comp.neg}, \ref{GB.neg},
\ref{noY.neg} and \ref{noT.neg} below then follow by $C^2$
approximation to $\Gamma$, via the argument on pp. 351--352 of
\cite{CG2}.

We shall compare $\Sigma$ with the {\em geodesic cone}
$C = C_p(\Gamma)$, which is formed from the minimizing geodesics
joining $p$ to points of $\Gamma$.  $C$ may naturally be given the
Riemannian metric $ds^2$ induced from $M^n$.  However, it should
be observed that $C$ with the metric $ds^2$ is not likely to be
relevant to the strongly stationary surface $\Sigma$.  In fact,
$\Sigma$ and the cone $C$ over its boundary inhabit different
regions of $M^n$, whose geometries are not related except by an
upper bound on curvatures, so that one should not expect any
useful comparison between them.  For these reasons, we shall endow
$C$ with a second metric $d\widehat s^2$ of constant Gauss
curvature $-\kappa^2$, such that the unit-speed geodesics from $p$
to points of $\Gamma$, which generate $C = C_p(\Gamma)$, remain
unit-speed geodesics in the metric $d\widehat s^2,$ and so that
$d\widehat s^2$ agrees with $ds^2$ at points of $\Gamma$
\cite{CG2}.  For clarity, we shall refer to the cone with this
hyperbolic metric as $\ch = \ch_p(\Gamma).$

More precisely, let $a_j, \ 1\leq j\leq m,$ be the smooth arcs  
of $\Gamma$, and let $A_j = C_p(a_j), \ 1\leq j\leq m,$
be the two-dimensional fans of $C_p(\Gamma)$.  On each $A_j$,  
let $\theta$ be a coordinate which is constant along each of the
radial geodesics through $p$, and such that
$\rho = \dist(\cdot, p)$ and $\theta$
form a local system of coordinates.  We have assumed that $\Gamma$
is nowhere tangent to the radial geodesic, which implies that
$\theta$ may be used as a regular parameter along the arc $a_j$.
Write $\rho =: r(\theta)$ for the corresponding values of
$\rho := \dist_M(p,\cdot)$ along $a_j$, and let $r(\theta)$ be  
extended to $C_p(\Gamma)$ so that it is constant along each  
radial geodesic.  
Then $\rho < r(\theta)$ elsewhere on $A_j$.  Note that under our
assumption, there holds $|dr/d\theta| < ds/d\theta$ along $\Gamma$.  
We may now write the metric $d\widehat s^2$ on $A_j$ as
$$
d\widehat s^2 = d\rho^2 +
\left[\left(\frac{ds}{d\theta}\Big |_\Gamma\right)^2 -
\left(\frac{dr(\theta)}{d\theta}\right)^2 \right]
\frac{\sinh^2\kappa\rho}{\sinh^2\kappa r(\theta)}\, d\theta^2.
$$

We may observe that,
along any radial geodesic, we have
$d\widehat s^2 = d\rho^2 = ds^2$.  In particular,
 if arcs $a_j$ and $a_k$ of $\Gamma$ share a
common endpoint $q$, then the hyperbolic metrics $d\widehat s^2$
defined on the fan $A_j$ and $d\widehat s^2$ defined on $A_k$
agree along their common edge, which is the minimizing geodesic
from $p$ to $q$.  That is, $d\widehat s^2$ makes
$\ch$ into a Riemannian polyhedron.

\noindent
\begin{thm}\label{dens.comp.neg}
Given a strongly stationary surface $\Sigma$ in $M^n$ of class
${\mathcal S}_{\Gamma}$, and a point $p$ of
$\Sigma \backslash \Gamma$, the following inequality holds:
\[
\Theta_{\Sigma}(p) \leq \Theta_{\ch_p(\Gamma)} (p).
\]
Moreover, equality implies that $\Sigma$ is a cone with totally
geodesic faces of constant Gauss curvature $-\kappa^2$.
\end{thm}

\pf
The proof is similar to the proof of Theorem
\ref{dens.comp}, with certain modifications.  The test
function $G(x)$ is taken to be $\log\tanh(\kappa\rho(x)/2),$
rather than $\log \rho(x).$  Since the faces $A_j$ of
$\ch_p(\Gamma)$ are locally isometric to the hyperbolic plane of
constant Gauss curvature $-\kappa^2$, with $\rho(x)$ corresponding
to the hyperbolic distance from a point, we may readily verify
that $G(x)$ is harmonic on the faces of $\ch_p(\Gamma)$ away from
$p$.  It follows from the trace formula \eqref{trace} and the
Hessian comparison theorem (p.~4 of \cite{SY}) that $G(x)$ is
subharmonic on the faces of $\Sigma$.  The factor $\frac{1}{\rho}$
appearing in boundary integrals in the proof of Theorem
\ref{dens.comp} is replaced by
$\frac{\kappa}{\sinh(\kappa\rho)}$, which is the derivative of $G$
with respect to $\rho$.  Note that
$-\kappa\,{\rm Length}(\ch\cap\partial B_{\varepsilon}(p))
/\sinh(\kappa\varepsilon)$
is equal to $- 2 \pi \Theta_{\ch}(p)$, independent of
sufficiently small $\varepsilon >0$.  If $e_n
=\overline\nabla\rho$ and $e_1,\dots,e_{n-1}$ form an orthonormal
frame on $M^n \backslash \{p\}$, then by the Hessian comparison
theorem $\overline\nabla^2_{e_i,e_i}G \geq \frac{\kappa^2 \cosh
\kappa\rho}{\sinh^2\kappa\rho}$ for $i = 1,\dots,n-1$, and
$\overline\nabla^2_{e_n,e_n}G =
-\frac{\kappa^2 \cosh \kappa\rho}{\sinh^2\kappa\rho}$
(See \cite{CG2}).                        
The remainder of the proof is as
in the proof of Theorem \ref{dens.comp}.
\qed

\noindent
\begin{thm} \label{GB.neg}
Let $\Gamma$ be a graph in $M^n$, and choose $p \in M^n$.  Then
the cone $\ch = \ch_p(\Gamma)$, with the hyperbolic metric
$d\widehat s^2$, satisfies the density estimate
\[
2\pi\Theta_\ch (p) \leq
- \sum_{k=1}^n \int_{a_k} \vec{k} \cdot \nu_C \,ds -
\kappa^2 \area \left(C_p(\Gamma)\right) +
\sum_k \sum_j \left(\frac{\pi}{2} - \beta_k^j\right),
\]
where $\nu_C$ is the outward unit normal vector to $C_p(\Gamma)$;
and at a vertex $q_j$ of $\Gamma$, $\beta_k^j$ is the angle
between the edge $a_k$ of $\Gamma$ and the minimizing geodesic
from $q_j \in \partial a_k$ to $p$.
\end{thm}

\pf
The proof is similar to the proof of Theorem \ref{GB} above.
We apply the Gauss-Bonnet formula \eqref{GBform} to the hyperbolic
cone $\ch = \ch_p(\Gamma)$, and find
$$
\int_{\ch\backslash B_\varepsilon(p)} K_\ch \, dA_\ch +
\int_{\ch\cap\partial B_\varepsilon(p)} \kh \, d\widehat s+
\int_{\Gamma_{\rm reg}} \kh \, d\widehat s +
\sum_k \sum_j \left(\frac{\pi}{2}-\widehat\beta_j^k \right) = 0,
$$
where $K_\ch \equiv -\kappa^2$ is the Gauss curvature of the faces
of $\ch$;  $\kh$ is the inward geodesic curvature along
$\partial\left(\ch\backslash B_\varepsilon(p)\right)$;  and
$\widehat\beta_j^k$ is the angle formed by the edge $a_k$ of
$\Gamma$ and the geodesic edge joining $p$ to
$q_j \in \partial a_k$, in the metric $d\widehat s^2$.  But along
$\partial B_\varepsilon(p) \cap \ch$, we have
$\kh \equiv -\kappa \coth \kappa\varepsilon$ by a standard computation
in the hyperbolic plane.  Along $\Gamma$, $d\widehat s^2 = ds^2$, so
that $\widehat\beta_j^k = \beta_j^k$.  Further, for each
$q \in \Gamma$, there holds $\kh(q) \leq k(q)$, the geodesic
curvature of $\Gamma$ in the cone $C_p(\Gamma)$ with the induced
metric $ds^2$ (see Proposition 4 of \cite{CG2}).  Thus
\begin{equation}\label{GBch.neg}
\kappa\coth\kappa\varepsilon\len(\partial B_\varepsilon(p)\cap\ch) \leq
\end{equation}
$$
-\kappa^2 \area(\ch \backslash B_\varepsilon(p))
+\int_{\Gamma_{\rm reg}} k \, ds +
\sum_k \sum_j \left(\frac{\pi}{2} - \beta_j^k \right).$$
Taking the limit as $\varepsilon \rightarrow 0$, we find
$$
2\pi \Theta_\ch(p) \leq
- \int_{\Gamma_{\rm reg}} \nu_C \cdot \vec{k} \, ds +
\sum_k \sum_j \left(\frac{\pi}{2} - \beta_j^k \right) -
\kappa^2 \area(\ch),
$$
since for all $q \in \Gamma$, $k(q) = -\nu_C \cdot \vec{k}(q).$
Finally, $\area(\ch) \geq \area(C),$ as may be proved by applying
Proposition 5 of \cite{CG2} to each face $A_k$ of $C$.
\qed

\bigskip

In order to state the following corollary and the next two theorems,
it will be useful to make the following

\noindent
\begin{defn}
$\Acal(\Gamma)$ is the {\em minimum cone area} of $\Gamma$:
$$
\Acal(\Gamma):= \min_{p\in \chull(\Gamma)} \area(C_p(\Gamma)).
$$
Here, the {\em convex hull} $\chull(\Gamma)$ of $\Gamma$ in $M$ is
the intersection of closed, locally geodesically convex subsets of
$M^n$ which contain $\Gamma$.
\end{defn}

\noindent
\begin{cor} \label{tc.bound.neg}
For a strongly stationary surface $\Sigma$ in a manifold $M^n$ with
sectional curvatures $K_M \leq -\kappa^2$,
the area-density estimate holds:
\[
2 \pi \Theta_{\Sigma}(p) \leq \tc(\Gamma) - \kappa^2 \Acal(\Gamma).
\]
Moreover, equality may only hold when $\Sigma$ is itself a cone
over $p$ with totally geodesic faces of constant Gauss curvature
$-\kappa^2$.
\end{cor}

\pf
Recall that Theorem \ref{GB.neg} estimates the hyperbolic cone density:
\begin{equation}\label{cdineq}
2\pi\Theta_\ch (p) \leq
- \sum_{k=1}^n \int_{a_k} \vec{k} \cdot \nu_C \,ds +
\sum_k \sum_j \left(\frac{\pi}{2} - \beta_k^j\right) -
\kappa^2 \area \Big(C_p(\Gamma)\Big).
\end{equation}
Since $\Sigma$ must lie in the convex hull $\chull(\Gamma)$ by the
maximum principle,  we have
$\area \left(C_p(\Gamma)\right) \geq \Acal(\Gamma)$.  Also,  
$|\int_{\Gamma_{\rm reg}} \vec{k} \cdot \nu_C \,ds| +
\sum_k \sum_j (\frac{\pi}{2} - \beta_j^k) \leq \tc(\Gamma).$
Therefore, the right-hand side of inequality
\eqref{cdineq} is $\leq \tc(\Gamma) - \kappa^2 \Acal(\Gamma)$,  
while according to Theorem \ref{dens.comp.neg}, the left-hand
side is $\geq 2 \pi \Theta_{\Sigma}(p)$.  Moreover, if equality
holds, then we must have equality in the conclusion of Theorem
\ref{dens.comp.neg}, implying that $\Sigma$ must be a cone
over $p$ with totally geodesic faces of constant Gauss curvature
$-\kappa^2$.
\qed

\medskip

In the following two theorems, the total curvature of $\Gamma$ is
``corrected" by subtracting $\kappa^2 \Acal(\Gamma)$.
Without this improved hypothesis, Theorems \ref{noY.neg} and
\ref{noT.neg} would have only extremely limited application for
$\Gamma$ of large diameter in manifolds $M^n$ of uniformly
negative sectional curvature
(see Example 2 of \cite{CG2}).

\noindent
\begin{thm}\label{noY.neg}
Suppose $\Gamma$ is a graph in $M^n$ with
$\tc(\Gamma) - \kappa^2 \Acal(\Gamma)\leq 3\pi$, and let
$\Sigma$ be a strongly stationary surface relative to $\Gamma$ in
the class ${\mathcal S}_{\Gamma}$.  Then $\Sigma$ is either an
embedded minimal surface;  or, a subset of a singular minimal cone
with an interior edge where three totally geodesic faces, of
constant Gauss curvature $-\kappa^2$, meet at equal angles.
\end{thm}

\pf
Given $p \in \Sigma$, Corollary \ref{tc.bound.neg}
above implies that
\[
2 \pi \Theta_{\Sigma}(p) \leq \tc(\Gamma) - \kappa^2 \Acal(\Gamma).
\]
Thus, the present hypothesis implies that
$\Theta_{\Sigma}(p) \leq \frac{3}{2}$, and that equality may only
hold when $\Sigma$ is a geodesic cone over $p$
and $\Sigma$ has totally geodesic faces of Gaussian curvature
$-\kappa^2$ (see Corollary \ref{tc.bound.neg}).
If
$\Theta_\Sigma (p) < 3/2$, then $\Sigma$ is embedded near $p$.  If
$\Theta_\Sigma (p) = 3/2$, then
$\Sigma$ is a geodesic cone, with tangent cone at $p$ congruent to the
Y stationary cone, and its faces are totally geodesic with Gauss
curvature $\equiv -\kappa^2$.  Since $\Sigma$ is a totally
geodesic cone of class ${\mathcal S}_{\Gamma}$, it is the
exponential image of its tangent cone at $p$.  It follows that the
exponential map of $M$ at $p$ maps a subset of the Y cone in $T_p M$ onto
$\Sigma$.
\qed

\noindent
\begin{thm}\label{noT.neg}
Suppose $\Gamma$ is a graph in $M^3$ with
$\tc(\Gamma) - \kappa^2 \Acal(\Gamma)\leq 2 \pi C_T$, and let
$\Sigma$ be an element of the regularity class
${\mathcal S}_{\Gamma}$, which is an
$(M, \varepsilon, \delta)$-minimal set   
with $\Gamma$ as its variational boundary.  Then $\Sigma$ is a
surface with possibly Y singularities but no other singularities
$p$, unless it is a geodesic cone over $p$ with totally geodesic
faces of constant Gauss curvature $-\kappa^2$, and having tangent
cone at $p$ equal to the T stationary cone.  
\end{thm}

\pf
Choose a point $p \in \Sigma$.  Then with respect to a  local
geodesic coordinate chart centered at $p$, the surface $\Sigma$ is an
$(M, \varepsilon, \delta)$-minimal set with $\varepsilon(r) = C
r^{\alpha}$ for some $C > 0$ and $\alpha >0$. Here we again apply
the set of results~\cite{T}(II.2 and II.3) to conclude that the
tangent cone $T_p \Sigma \subset T_p M^3 \cong \Re^3$ is area
minimizing and that the tangent cone can only be the plane, the
Y-cone or the T-cone.

As in the proof of Theorem
\ref{noY.neg}, we apply Corollary \ref{tc.bound.neg} to show that
either $\Theta_{\Sigma}(p) < C_T$;  or that $\Theta_{\Sigma}(p) = C_T$,
and $\Sigma$ is a geodesic cone over $p$ with
totally geodesic faces of constant Gauss curvature $-\kappa^2$,
which is  the image under the exponential map of $M$ at $p$ of the
T-cone.  If $\Theta_{\Sigma}(p) < C_T$, then the tangent cone to
$\Sigma$ at $p$ is either a plane or the Y stationary cone.  If
$T_p\Sigma$ is a plane, then $\Sigma$ is an embedded surface in a
neighborhood of $p$.  If $T_p\Sigma$ is the Y stationary cone,
then there are Y-type singularities along a curve passing through
$p$.  
\qed

\noindent
\begin{ram}
In Theorems \ref{noY.neg} and \ref{noT.neg}, the minimum cone area
$\Acal(\Gamma)$  may be replaced by
$$
\inf_{\chull(\Gamma)} \area(\ch_p(\Gamma)),
$$
which may be larger (and thus better).  See the proof of Theorem
\ref{GB.pos} below.  We have chosen to write Theorems \ref{noY.neg}
and \ref{noT.neg} in terms of the minimum cone area
$\Acal(\Gamma)$, since this quantity seems more closely
related to the geometry of $M$.    
(If $M$ has constant sectional curvature $-\kappa^2$, they are
equal.)
\end{ram}

\subsection{Ambient Curvature with Positive Upper Bound}\label{pos}

Throughout this subsection, weshall assume that $M^n$ is a
strongly convex Riemannian manifold whose sectional curvatures are
bounded above by a positive constant $\kappa^2$.  Consider a graph
$\Gamma \subset M^n$ and a surface $\Sigma$ of the regularity
class ${\mathcal S}_{\Gamma}$ which is strongly stationary with
respect to $\Gamma$.

Choose a point $p$ of $\Sigma$.  As in subsection \ref{neg}, we
shall assume that $\Gamma$ is nowhere tangent to the minimizing
geodesic from $p$.  The general cases of the results of this
subsection follow by $C^2$ approximation to $\Gamma$.

Since $M^n$ is strongly convex, the unique minimizing geodesic
joining $p$ to $q$ varies smoothly as a function of $q$.
Therefore, the geodesic cone $C = C_p(\Gamma)$, with the
Riemannian metric $ds^2$ induced from $M$, is a Riemannian
polyhedron enjoying the same smoothness as $\Gamma$.  This cone
will be given a second Riemannian metric $d\widehat s^2$, the
spherical metric, so that the faces of the cone have constant
Gauss curvature $\kappa^2$, so that the ambient distance $\rho$ to
the point $p$ remains equal to the distance in either metric
$ds^2$ or $d\widehat s^2$, and so that at points of $\Gamma$,
$d\widehat s^2 = ds^2$.  We may describe the spherical metric
at a point $q$ of $C$ as
$$
d\widehat s^2 = d\rho^2 +
\frac{\sin^2\kappa\rho}{\sin^2\kappa r(q)}
\left[\Big(ds\Big |_\Gamma\Big)^2 - \Big(dr(q)\Big)^2 \right].
$$
As in subsection \ref{neg}, $r(q)$ denotes $\rho(Q)$, the distance
in $M$ from $p$ to the point $Q$ of $\Gamma$ along the radial      
geodesic from $p$ passing through $q$;  also, the one-form        
$ds\Big |_\Gamma$   
has been extended to the cone so that it is invariant
under radial deformations   
Note that
$ds\Big |_\Gamma(\partial/\partial\rho) =
dr(\partial/\partial\rho) = 0$.  We use the notation
$\ch = \ch_p(\Gamma)$ for the cone $C$ with this spherical metric
$d\widehat s^2$.

In this section, it will be useful to state theorems in terms of a
{\em maximum} cone area, rather than the {\em minimum} cone area
which was of use in subsection \ref{neg}.  To account for the
positive sectional curvature which may occur in $M$, we will need
to add a term $\kappa^2 \Acalh(\Gamma)$ to the total curvature
$\tc(\Gamma)$.  The reader might object that, under certain
circumstances, such as when sectional curvatures comparable to
$\kappa^2$ appear only in a small part of $M^n$ and large parts
of the manifold $M$
actually have nonpositive sectional curvatures, this upper bound
may be much larger than the values which need to be considered in
Theorems \ref{noY.pos} and \ref{noT.pos} below.  However, when the
sectional curvatures of $M$ are nearly equal to the constant
$\kappa^2$, the theorems below are nearly sharp.

\noindent
\begin{defn}
$\Acalh(\Gamma)$ is the {\em maximum spherical cone area} of $\Gamma$:
$$
\Acalh(\Gamma)
:= \sup_{p\in \chull(\Gamma)}
\area(\ch_p(\Gamma)).
$$
\end{defn}

\noindent
\begin{thm}\label{dens.comp.pos}
Given a strongly stationary surface $\Sigma$ in $M^n$ of class
${\mathcal S}_{\Gamma}$, and a point $p$ of
$\Sigma \backslash \Gamma$, the following inequality holds:
\[
\Theta_{\Sigma}(p) \leq \Theta_{\ch_p(\Gamma)} (p).
\]
Moreover, equality implies that $\Sigma$ is a cone with totally
geodesic faces of constant Gauss curvature $\kappa^2$.
\end{thm}

\pf
Analogous to the proof of Theorem \ref{dens.comp.neg}, but using
$\log\tan(\kappa\rho(x)/2)$ as the test function $G(x)$ in place
of $\log\tanh(\kappa\rho(x)/2)$.
\qed

\noindent
\begin{thm} \label{GB.pos}
Let $\Gamma$ be a graph in $M^n$, and choose $p \in M^n$.  Then
the cone $\ch = \ch_p(\Gamma)$, with the spherical metric
$d\widehat s^2$, satisfies the density estimate
\[
2\pi\Theta_\ch (p) \leq
- \sum_{k=1}^n \int_{a_k} \vec{k} \cdot \nu_C \,ds +
\kappa^2 \area \left(\ch_p(\Gamma)\right) +
\sum_k \sum_j \left(\frac{\pi}{2} - \beta_k^j\right),
\]
where $\nu_C = \nu_\ch$ is the outward unit normal vector to
$C_p(\Gamma)$; and $\beta_k^j$ is the angle between the edge $a_k$
of $\Gamma$ and the minimizing geodesic in $M$ from
$q_j \in \partial a_k$ to $p$.
\end{thm}

\pf
The demonstration, which is based on the Gauss-Bonnet formula on
$\ch$, is highly analogous to the proof of Theorem \ref{GB.neg};
the {\em statement} has been modified, however, since in the
middle term on the right-hand side of equation \eqref{GBch.neg},
$\area(\ch)$ was multiplied by the non-positive $-\kappa^2$ and
could therefore be replaced in the conclusion of Theorem
\ref{GB.neg} with the smaller quantity $\area(C)$.  Here, however,
the Gauss curvature of $\ch$ is $\kappa^2$, which is positive, so
that the spherical area $\area(\ch)$ of the cone must remain on
the right-hand side of the inequality.
\qed

\noindent
\begin{cor} \label{tc.bound.pos}
The area density of a
strongly stationary surface $\Sigma$ in a manifold $M^n$ with
sectional curvatures $K_M \leq +\kappa^2$
satisfies the inequality:
\[
2 \pi \Theta_{\Sigma}(p) \leq \tc(\Gamma) + \kappa^2
\Acalh(\Gamma).
\]
Moreover, equality may only hold when $\Sigma$ is itself a cone
over $p$ with totally geodesic faces of constant Gauss curvature
$\kappa^2$.
\end{cor}

\pf
Theorem \ref{dens.comp.pos} estimates the density $\Theta_\Sigma(p)
\leq \Theta_{\ch_p(\Gamma)}(p)$.  Meanwhile, by Theorem \ref{GB.pos},  
\begin{eqnarray}
2\pi\Theta_\ch (p) &\leq&
- \sum_{k=1}^n \int_{a_k} \vec{k} \cdot \nu_C \,ds \\
&+&
\sum_k \sum_j \left(\frac{\pi}{2} - \beta_k^j\right) +
\kappa^2 \area \left(\ch_p(\Gamma)\right).
\end{eqnarray}
Since $\Sigma$ lies in the convex hull $\chull(\Gamma)$ by the
maximum principle, we have
$\area \left(\ch_p(\Gamma)\right) \leq \Acalh(\Gamma)$.  Also,
by definition of total curvature,
$|\int_{\Gamma_{\rm reg}} \vec{k} \cdot \nu_C \,ds| +
\sum_k \sum_j (\frac{\pi}{2} - \beta_j^k) \leq \tc(\Gamma).$
Therefore,
$2\pi\Theta_{\Sigma}(p)\leq\tc(\Gamma)+\kappa^2\Acal(\Gamma)$.
Moreover, if equality
holds, then we must have equality in the conclusion of Theorem
\ref{dens.comp.pos}, implying that $\Sigma$ must be a geodesic
cone over $p$ with totally geodesic faces of constant Gauss
curvature $+\kappa^2$.
\qed

The proofs of our final two theorems are completely analogous to
the proofs of Theorems \ref{noY.neg} and \ref{noT.neg}.

\noindent
\begin{thm}\label{noY.pos}
Suppose $\Gamma$ is a graph in $M^n$ with
$\tc(\Gamma) + \kappa^2 \Acalh(\Gamma)\leq 3\pi$, and let
$\Sigma$ be a strongly stationary surface relative to $\Gamma$ in
the class ${\mathcal S}_{\Gamma}$.  Then $\Sigma$ is either an
embedded minimal surface or a subset of a singular minimal cone
with an interior edge where three totally geodesic faces, of
constant Gauss curvature $\kappa^2$, meet at equal angles.
\end{thm}

\noindent
\begin{thm}\label{noT.pos}
Suppose $\Gamma$ is a graph in $M^3$ with
$\tc(\Gamma) + \kappa^2 \Acalh(\Gamma)\leq 2 \pi C_T$, and let
$\Sigma$ be a $(M, 0, \delta)$-minimal set with respect to $\Gamma$
in the regularity class ${\mathcal S}_{\Gamma}$.  Then $\Sigma$ is a
surface with possibly Y singularities but no other singularities
$p$, unless it is a geodesic cone over $p$ with totally geodesic
faces of constant Gauss curvature $\kappa^2$, and having tangent
cone at $p$ equal to the T stationary cone.
\end{thm}


\bigskip

\normalsize

\begin{tabbing}
aaaaaaaaaaaaaaaaaaaaaaaaaaaaaasssssssssssssssss \=
bbbbbbbbbbbbbbbbbbbbbbbbbbbbbbb\kill

Robert Gulliver\> Sumio Yamada\\
School of Mathematics \>Department of Mathematics\\
University of Minnesota \> University of Alabama, Birmingham\\
Minneapolis MN 55414 \> Birmingham, AL 35294\\
{\tt gulliver@math.umn.edu} \> {\tt yamada@math.uab.edu}\\
{\tt www.ima.umn.edu/\~{ }gulliver}\> {\tt www.math.uab.edu/yamada} \\
\end{tabbing}

\end{document}